\documentclass[12pt]{article}

\usepackage{color}
\usepackage{soul}
\usepackage{amsmath}
\usepackage{amssymb}
\usepackage{graphicx}
\usepackage{cite}
\usepackage{float}

\usepackage[ruled,linesnumbered]{algorithm2e}
\newtheorem{lemma}{Lemma}
\newtheorem{theorem}{Theorem}
\newtheorem{assumptions}{Assumption}
\newtheorem{rem}{Remark}
\newcommand{\fin}{\,\rule{1ex}{1.6ex}\,}

\date{ }

\title{Shanks and Anderson-type acceleration techniques for systems of nonlinear equations}
\author{Claude Brezinski\thanks{Universit\'{e} de Lille, CNRS, UMR 8524 - Laboratoire Paul Painlev\'{e},
F-59000 Lille, France. E--mail: {\tt Claude.Brezinski@univ-lille.fr}.}
\and
Stefano Cipolla\thanks{Universit\`a degli Studi di Padova,
	Dipartimento di Matematica ``Tullio Levi-Civita'',
	Via Trieste 63, 35121--Padova,
	Italy. E--mail: {\tt  stefano.cipolla87@gmail.com, michela.redivozaglia@unipd.it  }.}
\and Michela Redivo-Zaglia$^{\dag}$
\and
Yousef Saad\thanks{Dept. of Computer Science and Engineering, University of Minnesota, Mississippi National River and Recreation Area,
Minneapolis, MN 55455, USA.
E-mail: {\tt saad@cs.umn.edu}.}
}

\begin{document}

\maketitle

\begin{abstract}
This paper examines a number of extrapolation and acceleration methods,  and introduces a few
modifications of the standard Shanks transformation that deal with general sequences.
One of the goals of the paper is to lay out a general framework  that encompasses most of the known
acceleration strategies. The paper also considers the Anderson Acceleration method under a new light and exploits a
 connection with quasi-Newton methods, in order to establish  local linear convergence results of
a stabilized version of Anderson Acceleration method.
The methods are tested on a number of problems, including a few that arise from nonlinear Partial Differential Equations.

\vskip 2mm

\end{abstract}

\section{Introduction}
\label{sect1}

In numerical analysis and in applied mathematics, many applications lead to  sequences of numbers, vectors, matrices or
even tensors. When the sequence is slowly converging, or even diverging, and when one has only access to the
sequence and nothing else  (i.e., when it is produced by a ``black box"),
it is possible to transform it, by a {\it sequence transformation}, into a new sequence, which,
under some assumptions, converges faster than the original
one to the same limit. It was necessary to develop a variety of
such sequence transformations since, in fact,
it was  proved by Delahaye and Germain-Bonne \cite{bgbj} that a
universal sequence transformation able to accelerate all sequences, or
even all monotonically converging scalar ones, cannot exist. For a review,
see, for example, \cite{cbmrz,wimp,sidib,weni,dela,genesis,era}.

One way to transform a sequence into a faster converging one is to resort to {\it extrapolation}. Here, the transformation is built so that it
yields the exact limit of all sequences satisfying a certain algebraic relation. The set of these sequences is called the {\it kernel} of the transformation.
Among these, this paper focuses on
{\it Shanks transformation}   \cite{sha2} and a number of  its generalizations.
As we will see, this well-established method  transforms a sequence $(\mathbf s_n)$ into a set of sequences
$\{(\mathbf t_n^{(k)})\}$. Introduced by Shanks for scalar sequences
\cite{sha2}, it has been extensively studied, and extended, in various ways, to sequences of vectors, matrices, and tensors.
Here, we only consider the vector case.

All these extensions to vectors of the scalar Shanks transformation share the property that,
for a fixed value of $k$, $\mathbf  t_n^{(k)}=\mathbf s$ for all $n$
if the sequence $({\mathbf s}_n)$ of elements of $\mathbb R^{p}$ or $\mathbb C^{p}$ satisfies, for all $n$,
the following linear difference equation of order $k$
\begin{equation} \label{ker}
\alpha_0({\mathbf s}_n-{\mathbf s})+\cdots+\alpha_k({\mathbf s}_{n+k}-{\mathbf s})=0,
\end{equation}
where ${\mathbf s}$ is the limit of $({\mathbf s}_n)$ if it converges,
and is called its {\it antilimit} otherwise.  The numbers $\alpha_i$ are independent of $n$, and it is assumed that
$\alpha_0\alpha_k \neq 0$, so that the difference equation has the order $k$ exactly,
and $\alpha_0+\cdots+\alpha_k \neq 0$, so that $\mathbf s$ be uniquely defined. Thus, these conditions imply that $k$
cannot be replaced by a smaller value. It does not restrict the generality to assume that $\alpha_0+\cdots+\alpha_k=1$.
The set of sequences satisfying \eqref{ker} is called the {\it Shanks kernel}.
Among sequences in this kernel are those produced by the iterations $\mathbf s_{n+1}=M\mathbf s_n+\mathbf b$,
thus providing a link with Krylov subspace and Lanczos methods; see, in particular, \cite{birk,vec,sidibri,sidiprojection}.

\vskip 2mm

Besides their use in a number of different applications,
extrapolation techniques have recently been promoted as an   effective tool
also for problems related to the emerging field of Data
Science \cite{aspre,zhangboyd,cipollamulti,cipollashifted}. But
since there is often some confusion in the literature about the terminology used,
we would like clarify it -- using a high level of generality.
Specifically, we would like to draw a distinction between \textit{extrapolation methods},   \textit{sequence transformations},
and \textit{convergence acceleration methods}.
This distinction will help the reader  to better understand the
approaches described in Section \ref{sec:transformations_in_shanks} for building our sequence transformations.

Let $(\mathbf s_n)$ be a sequence of elements of a vector space $E$ on $\mathbb C$.
A common problem encountered in numerical analysis is to estimate the limit of this sequence from a certain number of its terms.
The problem can be solved  by an {\it extrapolation method} as follows \cite{cbth,ext}. Let
$$\boldsymbol \varphi: \mathbb N \times D \longmapsto E, \qquad D \subseteq \mathbb C^k,$$
be such that
$$\forall \mathbf b \in D, \quad \lim_{n \to \infty} \boldsymbol \varphi(n,\mathbf b)=0.$$
Let $V_{\boldsymbol \varphi}$ be the linear variety of sequences of elements of $E$ such that
$$\forall n, \quad \mathbf s_n=\mathbf s+\boldsymbol \varphi(n,\mathbf b),$$
with $\mathbf s \in E$. Obviously $\lim_{n \to \infty} \mathbf s_n=\mathbf s$.

By definition, if $(\mathbf s_n) \in V_{\boldsymbol \varphi}$, then,
$\forall n$, $\mathbf s=\mathbf s_n-\boldsymbol \varphi(n,\mathbf b)$.
Now, if $(\mathbf s_n) \notin V_{\boldsymbol \varphi}$, let us consider a sequence
$(\mathbf  t_n=\mathbf t+\boldsymbol \varphi(n,\boldsymbol \beta)) \in V_{\boldsymbol \varphi}$,
and impose that it satisfies
the {\it interpolation conditions} $\mathbf t_{n+i}=\mathbf s_{n+i}$ for $i=0,\ldots,k$.
The vector $\boldsymbol \beta\in D$ can be computed, assuming that it exists and is unique, in different ways as the solution
of a system of $k$ scalar equations that can be obtained as follows.
Let $E^*$ be the algebraic dual vector space of $E$, that is the vector space of linear functionals on $E$.
Let $\mathbf y,\mathbf y_1,\ldots,\mathbf y_k \in E^*$, and let $\langle \cdot, \cdot \rangle$ denote the duality product
between $E^*$ and $E$.
The first strategy consists in computing the vector $\boldsymbol \beta$ as the solution of the system
$$\langle \mathbf y_i,\mathbf s_{n+1}-\mathbf s_{n}\rangle=\langle \mathbf y_i, \boldsymbol \varphi(n+1,\boldsymbol \beta)-
\boldsymbol \varphi(n,\boldsymbol \beta)\rangle,
\quad i=1,\ldots,k.$$
In the particular case of Shanks transformation, writing this system in matrix form,
leads to a relation having the same structure as {\it Approach 3} in the {\it minimal residual approach}
of Section \ref{app3}, but with different indexes.

In the second strategy, the vector $\boldsymbol \beta$ is the solution of the system
$$\langle \mathbf y,\mathbf s_{n+i+1}-\mathbf s_{n+i}\rangle=\langle \mathbf y, \boldsymbol \varphi(n+i+1,\boldsymbol \beta)-
\boldsymbol \varphi(n+i,\boldsymbol \beta)\rangle,
\quad i=0,\ldots,k-1.$$
For Shanks transformation, this approach corresponds, in matrix form, to something similar to {\it Approach 6}
in the {\it topological approach} of Section \ref{app6}.

Then, in both cases, we set $\mathbf t=\mathbf s_n-\boldsymbol \varphi(n,\boldsymbol \beta)$. Since
$\mathbf t=\lim_{n \to \infty} \mathbf t_n$, it is an approximation of $\mathbf s$, and it has been obtained by {\it extrapolation}.
Obviously $\mathbf t$ depends on $n$ and $k$, and we will now denote it by $\mathbf t_n^{(p_k)}$
where $p_k+1$ denotes the number of elements of the initial sequence used in the process. Thus, when $n$ and $p_k$ vary,
the sequence $(\mathbf s_n)$ has been transformed into the set of sequences $\{(\mathbf t_n^{(p_k)})\}$. This procedure is named
an {\it extrapolation method}. An important remark to be made is that it is a purely algebraic procedure.
Richardson's and Romberg's methods, and Aitken's $\Delta^2$ process are such well known scalar extrapolation methods.
Thus, an extrapolation method results in a {\it sequence transformation} $T: (\mathbf s_n) \longmapsto (\mathbf t_n^{(p_k)})$
when either $p_k$ or $n$ is fixed, and the other index tends to infinity.
Conversely, most sequence transformations
can be interpreted as extrapolation methods. The variety $V_{\boldsymbol \varphi}$ is usually named the {\it kernel} of the
transformation $T$, and it is denoted ${\cal K}_T$.
If, when $n$ or $p_k$ tends to infinity, the sequence $(\mathbf t_n^{(p_k)})$ converges to $\mathbf s$ faster than the sequence
$(\mathbf s_n)$, the denomination {\it convergence acceleration method} is also used.
Let us mention that extrapolation
methods can also be applied to diverging sequences. They are often used for accelerating fixed point iterations, sometimes
coupled with a restarting strategy.

In this paper, instead of building Skanks transformation by computing the coefficients in \eqref{ker} as
the solution of a linear system in the usual way, we propose a new
optimization approach, based on minimization. This allows to easily introduce,
for sequences not belonging to the Shanks kernel, a unified framework that includes also regularized and preconditioned techniques.

Anderson Acceleration (AA) \cite{ander,ander1}, also called Anderson mixing, Pulay mixing or Direct Inversion
in the Iterative Subspace (DIIS) \cite{pul1}, in the computational physics and chemistry communities,
has been widely used and applied
to the solution of various fixed point problems over the last decades.
  The literature on this method is too broad to  allow for an exhaustive discussion
  but it suffices to search recent citations to this work to understand the truly exceptional renewed  interest in Anderson
  Acceleration   across many disciplines. A few of  the classical citations include the
papers  by Walker and Ni \cite{WalkerNi2011},  Higham and Strabi\'c \cite{a4},
Toth and Kelley \cite{toth2015convergence}, and by Fang and Saad \cite{fang2009two},
and  a few  papers that describe applications are
\cite{A,B,zhangboyd,kelley2018numerical,lupo2019convergence,C,pollock1}.

However, it is important here to stress
  that AA is not an extrapolation method in the exact sense defined above since it does not start from an
  arbitrary given sequence and transforms it into a new sequence. Instead it builds its own sequence step by step.
Anderson acceleration is in fact more akin to quasi-Newton techniques than to extrapolation. It was viewed as a form of
  secant method in the classic book by Ortega and Rheinboldt~\cite[pp. 204-205]{orte}.
  Its relations to secant type methods, specifically `multi-secant methods' was unraveled by Eyert
  \cite{eyert96}, and later exploited in  \cite{fang2009two} and also in  \cite{B}.
  In short, Anderson-Pulay mixing is a second order method whose goal is to
  accelerate a fixed point iteration.
If we were to allow the number of preview iterates used in the process to increase indefinitely we would
get something rather similar to a standard quasi-Newton method whose convergence would be superlinear at the limit.
This is not done in practice because of cost and numerical stability considerations.
However, a certain relation with the RRE method, which is an extrapolation method, exists, and AA can be recovered by
using the Coupled Shanks transformations, as explained in \cite{brzs}.
Due to this connection, we gave, in Section \ref{sec:anderson_type}, new
procedures in the style of Anderson acceleration, that are called
Anderson-Type Mixing (ATM in short).
Stabilized and regularized versions of AA will be also proposed.

\vskip 2mm

The outline of the paper is the following
\begin{itemize}
	\item[--]  In Section \ref{sec:transformations_in_shanks}, we
          present an overview of transformation techniques for sequences belonging to
          the Shanks kernel, and show how their limit or antilimit can be obtained exactly from
          these transformations. Four out of six of these techniques are presented in a new way that comes out from
          an optimization problem. Coupled sequences used in Section \ref{sec:anderson_type} are also described.

        \item[--] In Section \ref{sec:extr_techniques}, we present transformations based on the Shanks kernel. We show how to adapt and extend
 the idea proposed in \cite{aspre}  to our approaches. These modifications are specifically
          designed to accelerate general/nonlinear sequences which do not belong to the Shanks kernel.

	\item[--] In  Section \ref{sec:extrapolation_uses}, we  present the {\it Restarted} and the {\it
			Continuous-Updating} methods for exploiting the
		Shanks-based transformations presented in the previous section. In this way we are able to introduce
a unified  framework  able   to  encompass simultaneously  the newly introduced transformations and many of transformations already present
in the  literature.

\item[--] In Section \ref{sec:anderson_type}, we present new \textit{Anderson-Type Mixing}
  methods. We show how the classical AA fits into them. Then, we introduce preconditioning and a regularization strategies.
  Moreover, exploiting the connection with
  quasi-Newton methods, we prove the local linear convergence of a \textit{stabilized }version  of the classical AA, which allows us to
  substantiate theoretically the regularization strategy encompassed in the Anderson-type techniques previously presented in this section.

	\item[--] In Section \ref{sec:numerical_res}, we perform a
          comparative experimental study of some of the
          techniques proposed using, among other tests, a set of
          nonlinear problems arising from Partial Differential Equations (PDEs).
\end{itemize}

Let us explain our notation. Given a sequence $(\mathbf{s}_n)$, we set
$S_i^{(j)} = [{\mathbf s}_{i}, \ldots, {\mathbf s}_{i+j-1}] \in
\mathbb R^{p \times j}$.  Thus, the superscript $j$ corresponds to the
number of columns formed by the $p$-dimensional vectors of the
sequence $({\mathbf s}_{n})$, and the lower index $i$  is the index of the
first of these vectors in the sequence.  Whenever it is used, the forward difference operator $\Delta$ is applied to
the lower index, that is
$\Delta S_i^{(j)} =S_{i+1}^{(j)}-S_i^{(j)}=[\Delta {\mathbf s}_{i}, \ldots, \Delta
  {\mathbf s}_{i+j-1}]$, and similarly for $\Delta^2$. For a fixed
value of $k$, we denote by $\overline{S}_{i}^{(j)}$ the $kp \times j$
matrix formed by stacking the $k$ matrices
$S_i^{(j)},\ldots,S_{i+k-1}^{(j)}$ of dimension $p \times j$.  When
not explicitly indicated, the norm used is the Euclidean norm.
Throughout the paper, if not explicitly indicated, all matrices whose inverse is needed are assumed to be
nonsingular. If it is not the case, the pseudo-inverse may be used.

\section{Transformations for sequences in the Shanks kernel} \label{sec:transformations_in_shanks}
Let $(\mathbf s_n)$ be a sequence of vectors in $\mathbb R^p$ or $\mathbb C^p$ such that \eqref{ker} holds
 for a fixed value of $k$ and for all $n$.
Assuming, without loss of generality, that  $\sum_{i=0}^k\alpha_i=1$, then  we get from \eqref{ker}
\begin{equation} \label{eq:limiti_kernel}
	\alpha_0 \mathbf{s}_n+\cdots +\alpha_k \mathbf{s}_{n+k}= \mathbf{s}, \qquad \hbox{ for all } n \geq 0 .
\end{equation}
Alternatively, we can write
\begin{equation} \label{eq:lsq_equation}
\mathbf{s}_{n+k}- \sum_{j=0}^{k-1}\beta_j\Delta {\mathbf{s}}_{n+j}=\mathbf{s},
\end{equation}
with $\beta_j=\sum_{i=0}^j\alpha_i$ for $j=0,\ldots,k-1$ (note that the $\beta_i$'s
are defined in a slightly different way than in \cite[Sect. 3.1.3]{brzs}).

In Sections \ref{sec:shakns_RRE_equiv}  and \ref{sec:topological_approaches}, we show that when
$(\mathbf{s}_n)$ belongs to the Shanks kernel for a fixed value of $k$, it is possible to compute
exactly the limit or the antilimit of the sequence from a certain number $\ell_k$ (which depends on $k$ and on the transformation used)
of consecutive vectors of the sequence, where $\ell_k = k+2$ (for the Minimal residual approaches)  or
$\ell_k=2k+1$ (for the Topological approaches). For this purpose, we present
six different strategies for computing the coefficients
$\boldsymbol{\alpha}=(\alpha_0, \ldots, \alpha_k)^T$ or
$\boldsymbol{\beta}=(\beta_0, \ldots, \beta_{k-1})^T$. It should be reminded that $\boldsymbol{\alpha}$
and $\boldsymbol{\beta}$ are not dependent on $n$
if $(\mathbf s_n)$ satisfies \eqref{ker}, or \eqref{eq:limiti_kernel}, or \eqref{eq:lsq_equation}.
Four of these strategies (Approaches 1, 2, 4 and 5 below) are presented  as the solution of a
minimization problem. Approaches 1 and 4 proceed in what appears to be
a new way, not considered before in the literature
devoted to Shanks sequence transformations. Approaches 2 and 5 can be considered as particular cases
of the Least-squares strategy evoked in \cite[Sect. 3.1.3]{brzs}. These four strategies will be useful for the generalization presented in
Section \ref{sec:extr_techniques}.
Two of these strategies (Approaches 3 and 6 below) are already known since they enter into the framework of extrapolation methods as explained in Section \ref{sect1}, and are derived in Section \ref{sec:shakns_RRE_equiv}  and \ref{sec:topological_approaches}
by a purely algebraic process as
the solution of a linear system and they can be easily obtained by a modification of the Approaches 2 and 5. Moreover, as will be explained in Section \ref{sec:extr_techniques}, these
two strategies could also be included into the framework of the minimization by changing the metric of the norm. Approaches 3 and 6 will
be used in Section \ref{coupled}, where the notion of {\it coupled sequence}, defined in \cite{brzs}, is
invoked.

\vskip 0.3cm

Let us explain the idea behind the minimization used for finding the vector $\boldsymbol{\alpha}$ (since $\boldsymbol{\beta}$
is related to $\boldsymbol{\alpha}$, the idea is similar). This idea was introduced in \cite{aspre}, but it was not related to Shanks
transformations. In Section \ref{sec:extrapolation_uses} and the following ones,
our transformations are used to solve the fixed point problem $\mathbf s=G(\mathbf s)$ from iterates of the form
$\mathbf s_{n+1}=G(\mathbf s_n)$. Under some assumptions, it holds that
$\mathbf s_{n}-\mathbf s=(G'(\mathbf s))^n( \mathbf s_0-\mathbf s)+{\cal O}(\|\mathbf s_0-\mathbf s\|^2)$.
Thus, neglecting the terms of second order,
\[ \sum_{i=0}^k \alpha_i  \mathbf s_{n+i}-\mathbf s
  \approx \left(G'(\mathbf s) \right)^n\sum_{i=0}^k \alpha_i (G'(\mathbf s))^i( \mathbf s_0-\mathbf s)
  . \]
The idea is to minimize this error term. But
$\Delta \mathbf s_n \approx (G'(\mathbf s)-I)(\mathbf s_n-\mathbf s)$, and thus
\[ \sum_{i=0}^k \alpha_i \Delta \mathbf s_{n+i} \approx  (G'(\mathbf s)-I)
   (G'(\mathbf s))^n \sum_{i=0}^k \alpha_i (G'(\mathbf s))^i( \mathbf s_0-\mathbf s), \]
 which is similar to the expressions minimized for obtaining the vector
 $\boldsymbol{\alpha}$ in Approaches 1 and 4 below.

\vskip 2mm

When $\boldsymbol{\alpha}$ or
$\boldsymbol{\beta}$ has been computed, in any one of the ways described below,
 the vector $\mathbf{s}$ is  directly obtained by \eqref{eq:limiti_kernel} or  \eqref{eq:lsq_equation} as
 \begin{equation} \label{eq:shanks2}
\mathbf{s}=[ {\mathbf s}_{n+i},\ldots,  {\mathbf s}_{n+i+k}] \boldsymbol{\alpha}= S_{n+i}^{(k+1)}\boldsymbol{\alpha}, \qquad \mbox{for all} \;  i,
\end{equation}
or
\begin{equation} \label{eq:shanks2_SECOND}
\mathbf{s}= \mathbf{s}_{n+i+k}-[\Delta {\mathbf s}_{n+i},\ldots, \Delta {\mathbf s}_{n+i+k-1}]
\boldsymbol{\beta} = \mathbf{s}_{n+i+k}-\Delta {S}_{n+i}^{(k)}
\boldsymbol{\beta},\qquad \mbox{for all} \;  i.
\end{equation}

\begin{rem}
As can be seen, \eqref{eq:shanks2_SECOND} has the form of a Schur
complement
\begin{equation*}
\mathbf{u}= \mathbf{u}_{0}-[{\mathbf u}_{1},\ldots, {\mathbf u}_{k}]
A^{-1}\mathbf{v},
\end{equation*}
where $\mathbf{u}, \mathbf{u}_{0},{\mathbf u}_{1},\ldots, {\mathbf u}_{k}  \in  \mathbb{R}^{p}$,
$\mathbf{v} \in  \mathbb{R}^{k}$, and $A \in \mathbb{R}^{k \times k}$.
Several other expression in the sequel have the same form.

  From the extended Schur determinantal formula \cite{sch},
$\mathbf{u}$ can be expressed as a ratio of two determinants
\begin{equation*}
\mathbf{u}= \frac{
	\left|
	\begin{array}{ccc}
	{\mathbf u}_{0}& $\;$ & \mathbf{u}_{1} \cdots  \mathbf{u}_{k} \\
	\mathbf{v}  & $\;$  &   A
	\end{array}
	\right|}
{\left|
	A
	\right|} .
\end{equation*}
The determinant in the numerator is to be understood as the linear
  combination of the elements of its first row by applying the
  classical rules for expanding a determinant with
  respect its first row.
It is exactly through this connection that all the transformations given in \cite{brzs}
(Least-Squares strategy apart) have been defined.
\end{rem}

\subsection{Minimal residual approaches} \label{sec:shakns_RRE_equiv}
All the Minimal residual approaches described in this Section
for computing $\boldsymbol{\alpha}$ or $\boldsymbol{\beta}$  require
the knowledge of the $k+2$ vectors  $\mathbf{s}_n, \ldots, \mathbf{s}_{n+k+1}$.

\subsubsection{Approach 1}
\label{app1}

Writing \eqref{eq:limiti_kernel} for the indices $n$ and $n+1$ and subtracting, we obtain
\begin{equation*}
\alpha_0\Delta {\mathbf s}_n+\cdots+\alpha_k\Delta {\mathbf s}_{n+k}=0.
\end{equation*}
Then, one way to compute $\boldsymbol{\alpha}=(\alpha_0, \ldots, \alpha_k)^T$
 is to  solve the problem
\begin{equation} \label{eq:aspre_problem}
  \boldsymbol{\alpha} = \underset{\boldsymbol{\gamma} \in \mathbb{R}^{k+1}, {\mathbf e}^T \boldsymbol{\gamma}=1}{\arg \min}
  \| \Delta S_n^{(k+1)}\boldsymbol{\gamma} \|^2
\end{equation}
where
${\mathbf e}$ is the vector of all ones.
This is exactly the same relation introduced in \cite{aspre}, but obtained from a different starting point and without regularization.
The original paper by Pulay \cite{pul1} also solves the least squares problem with the same constraint that the sum of
the $\alpha_i$ equal to 1 by using Lagrange multipliers.

Observe that equation \eqref{eq:limiti_kernel} and the minimality of
$k$ ensure that $\dim \ker(\Delta S_n^{(k+1)})=1$. Hence, the solution of
problem \eqref{eq:aspre_problem} can  be obtained by normalizing
the unique vector in the kernel; alternatively, it can also be obtained as follows
(which leads to the SVD-MPE approach, see \cite{sidisvd})
\begin{equation} \label{eq:svd}
\boldsymbol{\alpha}=
\frac{\overline{\boldsymbol{\alpha}} }{{\mathbf e}^T\overline{\boldsymbol{\alpha}}}
\quad \hbox{where} \quad
  \overline{\boldsymbol{\alpha}} = \underset{\boldsymbol{\gamma} \in \mathbb{R}^{k+1}, \|\boldsymbol{\gamma}\|^2=1  }{\arg \min} 	\|
\Delta S_n^{(k+1)}\boldsymbol{\gamma} \|^2 .
\end{equation}

\subsubsection{Approach 2}
\label{app2}

Writing \eqref{eq:lsq_equation} for the indices $n$ and $n+1$ and subtracting,
we have

\begin{equation*}
	\Delta \mathbf{s}_{n+k}-\sum_{j=0}^{k-1}\beta_j [\Delta \mathbf{s}_{n+1+j}-\Delta \mathbf{s}_{n+j}]=\mathbf{0},
\end{equation*}
i.e., in compact form,

\begin{equation} \label{eq:RRE_strategy}
\Delta \mathbf{s}_{n+k}- \Delta^{2}S_n^{(k)}\boldsymbol \beta=\boldsymbol{0},
\end{equation}
where $\Delta^{2}S_n^{(k)}= [\Delta^2 \mathbf{s}_{n}, \ldots, \Delta^2 \mathbf{s}_{n+k-1}]$.

The vector $\boldsymbol{\beta}$ is solution of the problem
\begin{equation} \label{eq:RRE_solution}
\boldsymbol{\beta} =\underset{\boldsymbol{\eta} \in \mathbb{R}^k}{\arg \min} \|\Delta \mathbf{s}_{n+k}-
\Delta^{2}S_n^{(k)}\boldsymbol{\eta} ||^2,
\end{equation}
and therefore it can be obtained by solving the normal equations:
\begin{equation} \label{eq:RRE_solution_system}
({\Delta^{2}S_n^{(k)}})^T
\Delta^{2}S_n^{(k)} \boldsymbol{\beta} =
(\Delta^{2}S_n^{(k)})^T \Delta \mathbf{s}_{n+k},
\end{equation}
which leads to the strategy of the  Reduced Rank Extrapolation  (RRE)
due to Eddy \cite{rre} and Me\`sina \cite{rre1}.

\subsubsection{Approach 3}
\label{app3}

This approach generalizes the one seen in the preceding Section.
We consider a matrix $Y \in \mathbb{R}^{p \times k}$,
where $p$ is the dimension of the vectors of the sequence.
If we multiply \eqref{eq:RRE_strategy} by $Y^T$, it is  possible to
obtain the $\beta_i$ by solving the
following
 system that generalizes \eqref{eq:RRE_solution_system} which is obtained when
 $Y=\Delta^{2}S_n^{(k)}$
\begin{equation} \label{eq:tological_mmpe}
Y^T \Delta^{2}S_n^{(k)}\boldsymbol \beta= Y^T\Delta \mathbf{s}_{n+k},
\end{equation}
assuming that $\text{rank}(Y^T \Delta^{2}S_n^{(k)})=k$.

The best choice of the matrix $Y$ is a difficult problem which has not been studied yet. However,
some experimental results show that an appropriate choice of it can improve the convergence.
As shown, for example, in \cite{brzs}, particular choices of $Y$ yield several
existing extrapolation methods. Thus, the choice $Y=[{\mathbf  y}_1,\ldots,{\mathbf y}_k]$,
where the ${\mathbf y}_i$'s are $k$ linear independent vectors,
corresponds to the MMPE of Brezinski \cite{etopo} and Pugachev \cite{puga}
which can be recursively implemented by the $S\beta$-algorithm of Jbilou \cite{rrrr}.
The choice ${\mathbf   y}_i=\Delta {\mathbf s}_{n+i-1}$ leads to the MPE of Cabay and Jackson \cite{mpe}, and the RRE of  Me\`sina
\cite{rre1} and Eddy \cite{rre} is recovered with ${\mathbf   y}_i=\Delta^2 {\mathbf s}_{n+i-1}$ .

\subsection{Topological approaches} \label{sec:topological_approaches}

These approaches differ from those presented in Section \ref{sec:shakns_RRE_equiv} in that the algebraic equations
for computing the coefficients $\alpha_i$ or $\beta_i$ require more vectors of the sequence $({\mathbf s}_n)$, namely they
now need to utilize  the $2k+1$ vectors  $\mathbf{s}_n, \ldots, \mathbf{s}_{n+2k}$.

\subsubsection{Approach 4}
\label{app4}

Writing \eqref{eq:limiti_kernel} for the indices $n,\ldots, n+k$, and subtracting, we have
\begin{equation*}
	\alpha_0\Delta \mathbf{s}_{n+i}+\cdots +\alpha_k \Delta \mathbf{s}_{n+k+i}= 0, \qquad \hbox{ for } i=0,\ldots,k-1,
\end{equation*}
and the coefficients $\alpha_i$ are obtained by solving
\begin{equation}\label{eq:topological_no_normalization}
	\boldsymbol{\alpha}=\underset{\boldsymbol{\gamma} \in \mathbb{R}^{k+1}, {\mathbf e}^T\boldsymbol{\gamma}=1}{\arg \min} \|{\Delta \overline{S}_{n}^{(k+1)}} \boldsymbol{\gamma} \|^2
\end{equation}
where
\begin{equation*}
 \Delta \overline{S}_{n}^{(k+1)} = \left(
\begin{array}{cccc}
\Delta{\mathbf s}_n & \Delta{\mathbf s}_{n+1} & \cdots & \Delta{\mathbf s}_{n+k}\\
\Delta{\mathbf s}_{n+1} & \Delta{\mathbf s}_{n+2} & \cdots & \Delta{\mathbf s}_{n+k+1}\\
\vdots & \vdots & & \vdots \\
\Delta{\mathbf s}_{n+k-1} & \Delta{\mathbf s}_{n+k} & \cdots & \Delta{\mathbf s}_{n+2k-1}
\end{array}
\right) =
\left(
\begin{array}{c}
\Delta {S}_{n}^{(k+1)}\\
\Delta {S}_{n+1}^{(k+1)}\\
\vdots  \\
\Delta {S}_{n+k-1}^{(k+1)}
\end{array}
\right)
\in \mathbb{R}^{kp \times (k+1)}.
\end{equation*}

\subsubsection{Approach 5}
\label{app5}

The $\beta_i$'s can be computed by writing
\eqref{eq:lsq_equation} for the indices $n+k,\ldots, n+2k$, and subtracting. We have

\begin{equation*}
\Delta \mathbf{s}_{n+k+i}-\sum_{j=0}^{k-1}\beta_j \Delta^2 \mathbf{s}_{n+i+j}=\boldsymbol{0}, \qquad \hbox{ for } i=0,\ldots,k-1.
\end{equation*}
The coefficients $\beta_i$ are solution of the problem

\begin{equation} \label{eq:tea2_leqast_square}
	\boldsymbol{\beta}=\underset{\boldsymbol{\eta} \in \mathbb{R}^k}
{\arg \min}\|{\Delta \overline{S}_{n+k}^{(1)}}-{\Delta^2\overline{S}_{n}^{(k)}}\boldsymbol{\eta}\|^2
\end{equation}
where
\begin{equation*}
{\Delta \overline{S}_{n+k}^{(1)}}= \left(
\begin{array}{c}
\Delta \mathbf{s}_{n+k} \\
\vdots \\
\Delta \mathbf{s}_{n+2k-1}
\end{array}\right) \in
\mathbb{R}^{kp}\!, \;\; {\Delta^2\overline{S}_{n}^{(k)}}= \left(
\begin{array}{cccc}
\Delta^2{\mathbf s}_n & \Delta^2{\mathbf s}_{n+1} & \cdots & \Delta^2{\mathbf s}_{n+k-1}\\
\Delta^2 {\mathbf s}_{n+1} & \Delta^2{\mathbf s}_{n+2} & \cdots & \Delta^2{\mathbf s}_{n+k}\\
\vdots & \vdots & & \vdots \\
\Delta^2{\mathbf s}_{n+k-1} & \Delta^2{\mathbf s}_{n+k} & \cdots & \Delta^2{\mathbf s}_{n+2k-2}
\end{array}
\right) \in \mathbb{R}^{kp \times k},
\end{equation*}
that is
\begin{equation*}
\boldsymbol{\beta}=
(({\Delta^{2}\overline{S}_n^{(k)}})^T
\Delta^{2}\overline{S}_n^{(k)})^{-1}(\Delta^{2}\overline{S}_n^{(k)})^T \Delta {\overline{S}}_{n+k}^{(1)}.
\end{equation*}

\subsubsection{Approach 6}
\label{app6}

As in {Approach 3}, choosing $Y \in \mathbb{R}^{kp \times k}$,  we can
alternatively solve
\begin{equation}\label{eq:tea2}
Y^T{\Delta^2\overline{S}_{n}^{(k)}}\boldsymbol{\beta}=Y^T{\Delta \overline{S}_{n+k}^{(1)}}
\end{equation}
if $\textrm{rank}(Y^T{\Delta^2\overline{S}_{n}^{(2k-2)}})=k$.

When $Y=I_k \otimes \mathbf{y}$, for some $\mathbf{y} \in\mathbb{R}^{p}$, we recover the so called
{\it Topological Shanks transformation} that can be implemented recursively by the topological
$\varepsilon$-algorithms of Brezinski \cite{etopo} (in short TEA) or, more economically, by the
simplified topological $\varepsilon$-algorithms (in short STEA) \cite{simply,soft}.

\subsection{Coupled transformations}
\label{coupled}
We now recall the concept of {\it Coupled Sequences} introduced
  in \cite{brzs} since,  by using this extension, it is possible to link
  Anderson acceleration to the transformations based on the Shanks kernel.

  Given a sequence $(\mathbf{s}_n)$ belonging to the
Shanks kernel, a coupled sequence
$(\mathbf{c}_n)$ is a sequence which satisfies, for all $n$

\begin{equation*}
\alpha_0{\mathbf c}_n+\cdots+\alpha_k{\mathbf c}_{n+k}=0,
\end{equation*}
where the coefficients $\alpha_i$ are the same as in \eqref{eq:limiti_kernel}, or,
equivalentely
a sequence satisfying
  \begin{equation*}
  \mathbf{c}_{n+k}-\sum_{j=0}^{k-1}\beta_j \Delta \mathbf{c}_{n+j}=\mathbf{0}, \quad \hbox{ for all } n
  \end{equation*}
with the same coefficients $\beta_j$ as in \eqref{eq:lsq_equation}.
For example, the sequence ({$\mathbf{c}_n=\Delta^m \mathbf{s}_n$}) is a
sequence coupled to $(\mathbf{s}_n)$  for any $m \geq 1$.

By using a known coupled sequence, we can build additional generalizations of the Approaches 3 and 6,
which are recovered if we take ({$\mathbf{c}_n=\Delta \mathbf{s}_n$}),
and compute $\boldsymbol \beta$ as follows.
Let $C_n^{(k)} = [{\mathbf c}_{n}, \ldots, {\mathbf c}_{n+k-1}] \in \mathbb R^{p \times k}$.
Instead of \eqref{eq:tological_mmpe},   we solve the system
\begin{equation} \label{eq:coupled1}
Y^T \Delta C_n^{(k)}\boldsymbol{\beta}= Y^T \mathbf{c}_{n+k},
\end{equation}
where $Y \in \mathbb R^{p \times k}$.

Similarly, by defining the matrix $\overline{C}_{i}^{(j)}$ as made for $\overline{S}_{i}^{(j)}$, we can, instead of \eqref{eq:tea2}, solve
\begin{equation} \label{eq:coupled2}
Y^T{\Delta\overline{C}_{n}^{(k)}}\boldsymbol{\beta}=Y^T \overline{C}_{n+k}^{(1)},
\end{equation}
where now $Y \in \mathbb R^{kp \times k}$,

Particular choices of $Y$ and of the coupled sequence $(\mathbf{c}_n)$ give expressions similar to
those of well known methods (see \cite{brzs} for more details).

\section{Shanks-based transformations for general sequences} \label{sec:extr_techniques}

 The approaches described in the previous Section are all
 equivalent for a sequence belonging to the Shanks kernel
 and they yield the exact limit or antilimit.
 It is clear however, that this is an idealistic situation.
For extrapolating sequences that do not belong
to the Shanks kernel \eqref{ker}, we still write down the systems of
linear equations or the optimization problems giving the coefficients $\alpha_i$ or $\beta_i$
(which now depend of $k$ and $n$), and define a sequence transformation as the same linear combination
of terms as above.

In the sequel, for the extrapolated vector, we use the double indexing $\mathbf{t}_i^{(j)}$
that highlights the fact that the transformations,
require the $j+1$ elements
$\mathbf{s}_n, \ldots, \mathbf{s}_{n+j}$
of the sequence, in order to  compute $\boldsymbol{\alpha}$ or $\boldsymbol{\beta}$.

\begin{description}
\item[Minimal residual:]
In the case of the Minimal residual approaches there are $k+1$ vectors involved in the linear
combination. Thus, since to  compute $\boldsymbol{\alpha}$ or $\boldsymbol{\beta}$, we need
the $k+2$ vectors ${\mathbf s}_{n},\ldots,  {\mathbf s}_{n+k+1}$, we have only the
following
two different transformations, with the same $\boldsymbol{\alpha}$ and $\boldsymbol{\beta}$ (we denote the second transformation with a
tilde symbol over the ${\mathbf t}$)
\begin{itemize}
\item
$ \displaystyle
{\mathbf t}_n^{(k+1)}=
[{\mathbf s}_{n+1},\ldots,  {\mathbf s}_{n+k+1}] \boldsymbol{\alpha} =S_{n+1}^{(k+1)}\boldsymbol{\alpha}, \;
$ or equivalently \\
$ \displaystyle
{\mathbf t}_n^{(k+1)}=\mathbf{s}_{n+k+1}- [\Delta {\mathbf s}_{n+1},\ldots, \Delta {\mathbf s}_{n+k}]\boldsymbol{\beta}=
\mathbf{s}_{n+k+1}- \Delta S_{n+1}^{(k)}\boldsymbol{\beta}
$
\item
$ \displaystyle
{\widetilde {\mathbf t}}_n^{(k+1)}=
[{\mathbf s}_{n},\ldots,  {\mathbf s}_{n+k}] \boldsymbol{\alpha} =S_{n}^{(k+1)}\boldsymbol{\alpha}
$,
or equivalently\\
$ \displaystyle
{\widetilde{\mathbf t}}_n^{(k+1)}=\mathbf{s}_{n+k}- [\Delta {\mathbf s}_{n},\ldots, \Delta {\mathbf s}_{n+k-1}]\boldsymbol{\beta}=
\mathbf{s}_{n+k}- \Delta S_{n}^{(k)}\boldsymbol{\beta}
$
\end{itemize}
where   $\boldsymbol{\alpha} \in \mathbb{R}^{k+1}$ solves \eqref{eq:aspre_problem} or \eqref{eq:svd} (Approach 1) and
$\boldsymbol{\beta} \in \mathbb{R}^k$ solves \eqref{eq:RRE_solution} or \eqref{eq:tological_mmpe} (Approaches 2 or 3), or \eqref{eq:coupled1} (coupled approach).
\item[Topological:]
In the Topological case, there are again $k+1$ vectors involved in the linear
combination, but since we need
the $2k+1$ vectors ${\mathbf s}_{n},\ldots,  {\mathbf s}_{n+2k}$
to compute $\boldsymbol{\alpha}$ or $\boldsymbol{\beta}$,  we have $k+1$
different transformations (depending on the choice of the vectors used in the linear combination),
but with  the same $\boldsymbol{\alpha}$ and $\boldsymbol{\beta}$, and  we have
\begin{itemize}
\item
$ \displaystyle
	\mathbf{t}_{n,i}^{(2k)}=[{\mathbf s}_{n+i},\ldots,  {\mathbf s}_{n+i+k}] \boldsymbol{\alpha} =S_{n+i}^{(k+1)}\boldsymbol{\alpha}, \qquad
\mbox{for} \; i=0, \ldots, k$
\\
or equivalently\\
$ \displaystyle
	{\mathbf t}_{n,i}^{(2k)}=\mathbf{s}_{n+i+k}-
[\Delta {\mathbf s}_{n+i},\ldots, \Delta {\mathbf s}_{n+i+k-1}] \boldsymbol{\beta}=
\mathbf{s}_{n+i+k}-\Delta S_{n+i}^{(k)}\boldsymbol{\beta}, \qquad
\mbox{for} \; i=0, \ldots, k$
\end{itemize}
where $\boldsymbol{\alpha} \in \mathbb{R}^{k+1}$ solves \eqref{eq:topological_no_normalization} (Approach 4) and  $\boldsymbol{\beta} \in \mathbb{R}^k $
given by \eqref{eq:tea2_leqast_square} or \eqref{eq:tea2} (Approaches 5 or 6),
or \eqref{eq:coupled2} (coupled approach).
\\
Among all the possible linear combinations, it seems more appropriate
to use those involving the last available vector of the sequence, that is the
transformation with $i=k$ that uses
$\mathbf{s}_{n+2k}$. In the sequel, for simplifying the notation we will set
$\mathbf{t}_{n}^{(2k)}=\mathbf{t}_{n,k}^{(2k)}$.
\end{description}

\vskip 0.3cm

Of course, if the sequence belongs to the Shanks kernel, all the preceding transformations are equivalent and give the same result, that is $\mathbf s$.

\vskip 0.3cm

Now, let us show how to adapt and extend to our approaches
 the idea proposed in \cite{aspre}.
 All the transformations summarized at the beginning of this Section can be used,
 and the only change deals with the computation of the coefficients $\alpha_i$ or the
 $\beta_i$.  In \cite{aspre},  in
order to overcome the problems due to the ill-conditioning of
problem \eqref{eq:aspre_problem} (our Approach 1)
the authors consider
the following regularized problem for the computation of the $\alpha_i$ in the minimal residual approach
\begin{equation*}
\boldsymbol{\alpha}_{\lambda} =  \underset{ \boldsymbol{\gamma} \in \mathbb{R}^{k+1}, {\mathbf e}^T\boldsymbol{\gamma}=1}{\arg \min}	
\left(\| \Delta S_n^{(k+1)}\boldsymbol{\gamma} \|^2+\lambda \|\boldsymbol{\gamma}\|^2 \right),
\end{equation*}
with $\lambda \in \mathbb{R}$, and
whose solution is (assuming that ${\Delta S_n^{(k+1)}}$ is of full rank)
\begin{equation*}
\boldsymbol{\alpha}_{\lambda} = \frac{(({\Delta S_n^{(k+1)}})^T \! \Delta S_n^{(k+1)}+\lambda I)^{-1}{\mathbf e}}
{{\mathbf e}^T (({\Delta S_n^{(k+1)}})^T \! \Delta S_n^{(k+1)}+\lambda I)^{-1}{\mathbf e}} .
\end{equation*}
Observe that an alternative approach would be to change the metric in the evaluation of the norm, i.e., instead of using
the Euclidean norm,   solve the problem
\begin{equation} \label{eq:aspre_problem_prec}
\boldsymbol{\alpha}_{M,\lambda} =  \underset{\boldsymbol{\gamma} \in \mathbb{R}^{k+1}, {\mathbf e}^T\boldsymbol{\gamma}=1 }{\arg \min} 	
\left(\|
\Delta S_n^{(k+1)}\boldsymbol{\gamma} \|_M^2+\lambda \|\boldsymbol{\gamma}\|^2\right),
\end{equation}
where $\| \mathbf{x} \|_M^2 =  (\mathbf{x},M\mathbf{x})$
and $M \in \mathbb{R}^{p \times p}$ is a positive definite matrix.

  In what follows we will need  $M$ to be
  positive semi-definite only instead of positive definite. In this case
  $\| \cdot\|_M$ is a semi-norm but we abuse the terminology by calling it a
  `norm'.

With this, we have the following lemma.
\begin{lemma}
	The solution of problem \eqref{eq:aspre_problem_prec} is
	\begin{equation} \label{st0}
	\boldsymbol{\alpha}_{M,\lambda} =\frac{(({\Delta S_n^{(k+1)}})^T
\! M\Delta S_n^{(k+1)}+\lambda I)^{-1}{\mathbf e}}{{\mathbf e}^T (({\Delta S_n^{(k+1)}})^T \! M\Delta S_n^{(k+1)}+\lambda I)^{-1}{\mathbf e}},
	\end{equation}
	and the corresponding extrapolated vector is
	\begin{equation}  \label{st}
	{\mathbf t}_n^{(k+1)}=S_{n+1}^{(k+1)}\boldsymbol{\alpha}_{M,\lambda} \quad
{\it or} \quad {\widetilde{\mathbf t}}_n^{(k+1)}=S_{n}^{(k+1)}\boldsymbol{\alpha}_{M,\lambda} .
	\end{equation}
\end{lemma}
{\bf Proof:}
	The result follows by writing the problem \eqref{eq:aspre_problem_prec} as
	
	\begin{equation*}
		\boldsymbol{\alpha}_{M,\lambda} = \underset{\boldsymbol{\gamma} \in \mathbb{R}^{k+1}, {\mathbf e}^T\boldsymbol{\gamma}=1 }{\arg \min} 	\left(\boldsymbol{\gamma}^T(\Delta S_n^{(k+1)})^T \! M\Delta S_n^{(k+1)}\boldsymbol{\gamma} +\lambda \boldsymbol{\gamma}^T\boldsymbol{\gamma}
\right),
	\end{equation*}
	and by applying a technique analogous to that used in \cite{aspre}. From \eqref{eq:shanks2} we obtain \eqref{st}.
\fin

\vskip 0.2cm

Motivated by the equivalence of all the approaches described in
Section \ref{sec:shakns_RRE_equiv} for sequences in the Shanks kernel, we can thus
introduce the following problem
\begin{equation} \label{eq:RRE_problem_prec}
\boldsymbol{\beta}_{M,\lambda}=\underset{\boldsymbol{\eta} \in \mathbb{R}^k}{\arg \min}
\left(\|\Delta \mathbf{s}_{n+k}- \Delta^{2}S_n^{(k)}\boldsymbol{\eta} ||_M^2+\lambda \|\boldsymbol{\eta}\|^2 \right),
\end{equation}
where $M$ is a semi-positive definite matrix.
Referring to the gradient of the function $g(\boldsymbol{\eta})=\|\Delta \mathbf{s}_{n+k}- \Delta^{2}S_n^{(k)}
 \boldsymbol{\eta} ||_M^2+\lambda \|\boldsymbol{\eta}\|^2$, the solution of \eqref{eq:RRE_problem_prec} is given by
\begin{equation} \label{eq:beta_expression0}
\boldsymbol{\beta}_{M,\lambda}=(({\Delta^{2}S_n^{(k)}})^T \! M\Delta^{2}S_n^{(k)}+\lambda I)^{-1}(\Delta^{2}S_n^{(k)})^T \!
M \Delta \mathbf{s}_{n+k},
\end{equation}
and hence,  the corresponding extrapolated vector is
\begin{equation} \label{eq:beta_expression}
{\mathbf t}_n^{(k+1)}=\mathbf{s}_{n+k+1}- [\Delta {\mathbf s}_{n+1},\ldots, \Delta {\mathbf s}_{n+k}]\boldsymbol{\beta}_{M,\lambda},
\end{equation}
or
\begin{equation*}
\widetilde{\mathbf t}_n^{(k+1)}=\mathbf{s}_{n+k}- [\Delta {\mathbf s}_{n},\ldots, \Delta {\mathbf s}_{n+k-1}]\boldsymbol{\beta}_{M,\lambda},
\end{equation*}
In particular, if $M=YY^T$ where $Y\in \mathbb{R}^{ p\times k }$ is a given matrix  and $\lambda = 0$, we have,
\begin{equation*}
	({\Delta^{2}S_n^{(k)}})^TY \left(Y^T\Delta^{2}S_n^{(k)}
\boldsymbol{\beta}_{YY^T,0} -Y^T \Delta \mathbf{s}_{n+k}\right)=
{\mathbf 0.}
\end{equation*}

When  $\text{rank}(Y^T\Delta^{2}S_n^{(k)})=k$, we see that Approach 3
\eqref{eq:tological_mmpe} is a particular case of problem
\eqref{eq:RRE_problem_prec}.  As we already observed, different
choices of $Y \in \mathbb{R}^{p \times k}$ give rise to different {\it
  acceleration} performances for different type of sequences.

Similarly, following the idea of the topological approaches of Section \ref{sec:topological_approaches},
we consider the problems
\begin{equation}\label{eq:topological_no_normalization_regularized}
\boldsymbol{\alpha}_{M, \lambda }=\underset{\boldsymbol{\gamma} \in \mathbb{R}^{k+1}, {\mathbf e}^T\boldsymbol{\gamma}=1 }{\arg \min}
\left(\|{\Delta \overline{S}_{n}^{(k+1)}} \boldsymbol{\gamma} \|_M^2+\lambda\|\boldsymbol{\gamma}\|^2 \right),
\end{equation}
or
\begin{equation} \label{eq:TES_problem_prec}
\boldsymbol{\beta}_{M,\lambda}=\underset{\boldsymbol{\eta} \in \mathbb{R}^k}{\arg \min}
\left(\|{\Delta \overline{S}_{n+k}^{(1)}}-{\Delta^2\overline{S}_{n}^{(k)}}\boldsymbol{\eta}\|^2_M+\lambda \|\boldsymbol{\eta}\|^2\right),
\end{equation}
where, in both cases, $M  \in \mathbb{R}^{kp \times kp}$  is a semi-positive definite matrix.
The solution of \eqref{eq:topological_no_normalization_regularized} is
\begin{equation}\label{eq:solstar}
\boldsymbol{\alpha}_{M,\lambda} =\frac{(({\Delta \overline{S}_n^{(k+1)}})^T
\! M {\Delta \overline{S}_n^{(k+1)}}+\lambda I)^{-1}{\mathbf e}}{{\mathbf e}^T (({\Delta \overline{S}_n^{(k+1)}})^T \! M{\Delta \overline{S}_n^{(k+1)}})+\lambda I)^{-1}{\mathbf e}},
\end{equation}
and the corresponding extrapolated vector is
\begin{equation}\label{eq:topological_extrapolated}
\mathbf{t}_n^{(2k)}= S_{n+k}^{(k+1)} \boldsymbol{\alpha}_{M,\lambda}.
\end{equation}

The solution of problem \eqref{eq:TES_problem_prec} is
\begin{equation*}
\boldsymbol{\beta}_{M,\lambda}=(({\Delta^2\overline{S}_{n}^{(k)}})^T \! M{\Delta^2\overline{S}_{n}^{(k)}}+\lambda I)^{-1}({\Delta^2\overline{S}_{n}^{(k)}})^T
\! M {\Delta \overline{S}_{n+k}^{(1)}}
\end{equation*}
and the corresponding extrapolated vector is
\begin{equation}\label{eq:topological_extrapolated_beta_new}
{\mathbf t}_n^{(2k)}=\mathbf{s}_{n+2k}-[\Delta \mathbf{s}_{n+k},\ldots,\Delta \mathbf{s}_{n+2k-1}]\boldsymbol{\beta}_{M,\lambda}.
\end{equation}

We set $M=YY^T$ with
$Y \in \mathbb{R}^{kp \times k}$ and $\text{rank}(Y^T{\Delta^2\overline{S}_{n}^{(k)}})=k$. If
$\lambda = 0$, we see that the Approach 6 is a particular case of the problem \eqref{eq:TES_problem_prec}. If
$Y=I_{k}\otimes \mathbf{y}$, for some $\mathbf{y} \in\mathbb{R}^{p}$, we obtain
a method similar
to the Topological Shanks transformation
\cite{etopo}.

We refer the reader to Section \ref{sec:numerical_res} for a discussion of different possible strategies for the selection of the regularization parameter $\lambda$.

\section{Possible uses of acceleration strategies} \label{sec:extrapolation_uses}
In this Section and in the following one, we consider the solution of the fixed point problem $G(\mathbf{s})=\mathbf{s}$.
There are three ways to proceed.

The simplest way is to use an extrapolation method.
  The vectors $\mathbf{s}_n$ are generated one by one by Picard's iteration as
$\mathbf{s}_{n+1}=G(\mathbf{s}_n), n=0,1,\ldots$,  from a given $\mathbf s_0$. The extrapolation method is
applied
after each computation of a new vector $\mathbf{s}_n$ by using a certain number of the
preceding Picard's iterates to produce a completely new extrapolated sequence.
This procedure is called the \textit{Acceleration Method}
 but it is not used in this paper  (see \cite{soft} for details).

The second way consists in computing a certain number of Picard's iterates, then to use these
  in one of the extrapolation
  techniques introduced in Section \ref{sec:extr_techniques}, and finally to restart the Picard's
  iterates from the extrapolated vector that has been
 obtained. This is the {\it Restarted} method treated below.

 In the third way, the process builds its own sequence step by step. Each term of
 the sequence is obtained by combining, in a certain manner, Picard's iterates, preceding terms
 of the sequence and extrapolated ones.
 We will focus on three possible algorithms of this type that are termed {\it Continuous-Updating},
 presented in this Section, the {\it Anderson-type} and
 the {\it Periodic Anderson-type} methods, both discussed in Section \ref{sec:anderson_type}.
The difference between these procedures lies in the way in which previous  iterates are combined together in the process
  to obtain a new vector.

\subsection{Restarted method}
In this methodology,
already described, for example, in \cite{cb1,soft,gek}, a certain number of Picard's iterates
are produced, an extrapolation strategy is then applied to
them, and the Picard's iterates are {\it restarted} from the extrapolated vector; see Algorithm~\ref{RMtable}.
The sequence of the successive extrapolated terms
will be denoted by $(\mathbf{x}_j)$.

\vskip 2mm

\begin{algorithm}[H]
\caption{The Restarted Method (RM).} \label{RMtable}
	\DontPrintSemicolon
	\KwIn{Choose  $M$, $\lambda$, $k$, and ${\mathbf x}_0 \in \mathbb{R}^p$.}
	\For{$j=0,1,\ldots$}{Set ${\mathbf s}_0= {\mathbf x}_j$\;
			\For{$i=1,\ldots,\ell_k-1$ \it (basic or inner iterations)}{
			Compute ${\mathbf s}_i=G({\mathbf s}_{i-1})$ \;}
			Compute ${\mathbf t}_0^{(\ell_k-1)}$ using \eqref{st}~or~\eqref{eq:beta_expression}~or~\eqref{eq:topological_extrapolated}~or~\eqref{eq:topological_extrapolated_beta_new} \;
			Set ${\mathbf x}_{j+1}={\mathbf t}_0^{(\ell_k-1)}$\;
			}
\end{algorithm}

\vskip 2mm

Observe that $\ell_k=k+2$ if we use
\eqref{st} or \eqref{eq:beta_expression}, and $\ell_k=2k+1$ if we use
\eqref{eq:topological_extrapolated}
or \eqref{eq:topological_extrapolated_beta_new}.
In the particular case of \eqref{eq:beta_expression}, we have
\begin{equation*}
{\mathbf t}_0^{(k+1)}=\mathbf{s}_{k+1}-G_{M,\lambda}  \Delta \mathbf{s}_{k}
\end{equation*}
where
\begin{equation*}
G_{M,\lambda}=[\Delta {\mathbf s}_{1},\ldots, \Delta {\mathbf s}_{k}]
(({\Delta^{2}S_0^{(k)}})^T \!M\Delta^{2}S_0^{(k)}+\lambda
  I)^{-1}(\Delta^{2}S_0^{(k)})^T \!M,
\end{equation*}

Setting $\mathbf{f}_k= G(\mathbf{s}_k)-\mathbf{s}_k=\Delta \mathbf{s}_k$, we have
\begin{equation*}
{\mathbf t}_0^{(k+1)}=\mathbf{s}_{k+1}- G_{M,\lambda}  \mathbf{f}_k.
\end{equation*}
Therefore, we can interpret the Restarted Method as a {\it
  cyclic projection method} (see \cite{projproj} and \cite{proj}
for the linear case) for the solution of the problem $F(\mathbf{s})=0$
where $F(\mathbf{s})=G(\mathbf{s})-\mathbf{s}$.

The idea of the RM, that is to interleave a certain number of Picard's iterates with one extrapolation step, can also be used
in the {\it Continuous-Updating} and in the {\it Anderson-type} methods (see Section \ref{sec:anderson_type} where
a general `periodic' algorithm of this type is presented).

\vskip 2mm

A particular case of the RM is the {\it Generalized Steffensen Method} (GSM)
which corresponds to the case where the dimension of the projection space coincides with the dimension of the system, that is for $k=p$.
Under some assumptions,  when $\lambda =0$ and $M=I$,
the sequence $({\mathbf x}_j)$ obtained by the GSM
asymptotically converges quadratically to the
fixed point $\mathbf s^*$ of $G$ even if $G$  is not a contraction.
The GSM is a generalization of the well-known Steffensen method
\cite{stef} when $p=1$. It was first proposed by Brezinski \cite{cb1} and Gekeler \cite{gek}
for the case of the vector $\varepsilon$-algorithm, but there was a gap in their proofs as in that
of Skelboe for the MPE \cite{skel}  as noticed in \cite{sfs}
The first complete proof of the quadratic convergence of the GSM was given by Ortega and Rheinbolt  \cite[p. 373]{orte} for Henrici's method \cite[pp. 115 ff.]{henrici}
(a particular case of the MMPE),
Le Ferrand \cite{hlf} for the first Topological Shanks transformation of Brezinski \cite{etopo}, and Jbilou and Sadok for the
MPE and the RRE \cite{rrrr}.

\subsection{Continuous-Updating method} \label{sec:continuos_updating}

In this approach, the sequence is {\it continuously accelerated} by computing a new basic iterate at each step, using it in the
extrapolation process, and, after the computation of the extrapolated vector, replacing the new basic iterate computed before by it.
Thus, when compared with the original fixed point sequence,  the continuous updating scheme builds a completely new
sequence whose iterates replace those of the original sequence.

We start with the Minimal residual approach for computing the $\alpha_i$.
We have
the following  {\it Continuous-Updating Method}  (Algorithm \ref{AAtable_1})

\vskip 2mm
\begin{algorithm}[H]
	\caption{Continuous-Updating Method (CU) with $\boldsymbol{\alpha}_{M,\lambda}$.} \label{AAtable_1}
	\DontPrintSemicolon
	\KwIn{Choose $M$, $\lambda$,  $m \in \mathbb{N}, m \geq 1$,  ${\mathbf s}_0 \in \mathbb{R}^p$.}
	\For{$j=0,1,\ldots$}{
		Set $m_j=\min(m,j)$\;
        Compute ${\mathbf{s}}_{j+1}= G(\mathbf{s}_{j})$ {\it (Picard iteration)}\;
   Set $ {S}_{j-m_j}^{(m_j+1)}=[ \mathbf{s}_{j-m_j}, \ldots, \mathbf{s}_{j}]$\;
        Compute $\boldsymbol{\alpha}_{M,\lambda}$ using~\eqref{st0}~and~$\Delta {S}_{j-m_j}^{(m_j+1)}$\;
        Compute $\widetilde{\mathbf t}_{j-m_j}^{(m_j+1)}=
{S}_{j-m_j}^{(m_j+1)}\boldsymbol{\alpha}_{M,\lambda}$\;
		Set ${\mathbf s}_{j+1}= \widetilde{\mathbf t}_{j-m_j}^{(m_j+1)}$ \;
	}
\end{algorithm}

\vskip 2mm

Algorithm \ref{AAtable_11} listed next, uses formulas
(\ref{eq:beta_expression0}--\ref{eq:beta_expression}),
(i.e. the ${\beta_i}$'s, are computed by
\eqref{eq:beta_expression0}, that solve the problem \eqref{eq:RRE_problem_prec})

\vspace{2mm}

\begin{algorithm}[H]
	\caption{Continuous-Updating Method (CU) with $\boldsymbol{\beta}_{M, \lambda}$.} \label{AAtable_11}
	\DontPrintSemicolon
	\KwIn{Choose $M$, $\lambda$,  $m \in \mathbb{N}, m \geq 1$,  ${\mathbf s}_0 \in \mathbb{R}^p$.}
	Compute  ${\mathbf s}_1=G({\mathbf s}_0)$ \;
	\For{$j=1,2,\ldots$}{
		Set $m_j=\min(m,j)$\;
        Compute  ${\mathbf s}_{j+1}=G({\mathbf s}_j)$ {\it (Picard iteration)}\;
        Set $ \Delta {S}_{j-m_j}^{(m_j)}=[ \Delta \mathbf{s}_{j-m_j}, \ldots, \Delta \mathbf{s}_{j-1}]$\;
        Compute $\boldsymbol{\beta}_{M,\lambda}$ using~\eqref{eq:beta_expression0}~and~$\Delta^2 {S}_{j-m_j}^{(m_j)}$\;
        Compute $\widetilde{\mathbf t}_{j-m_j}^{(m_j+1)}=
\mathbf{s}_j - \Delta {S}_{j-m_j}^{(m_j)}\boldsymbol{\beta}_{M,\lambda}$\;
Set ${\mathbf s}_{j+1}= \widetilde{\mathbf t}_{j-m_j}^{(m_j+1)}$ \;
	}
\end{algorithm}

\vskip 2mm

As in the preceding algorithm, the new fixed point iterate ${\mathbf
  s}_{j+1}$ is used only for computing $\boldsymbol{\beta}_{M,
  \lambda}$. Thereafter, this iterate is not used in the linear combination for computing the
extrapolated vector as it is replaced by the extrapolated one that is computed.
\vskip 2mm

It is possible to highlight the connection between acceleration techniques and the projection framework. We define
\begin{equation*}
	G^{(j)}_{M,\lambda}=[\Delta {\mathbf s}_{j-m_j},\ldots, \Delta {\mathbf s}_{j-1}](({\Delta}^2S_{j-m_j}^{(m_j)})^T \! M{\Delta}^2S_{j-m_j}^{(m_j)}+\lambda I)^{-1}({\Delta}^2S_{j-m_j}^{(m_j)})^{T} \!M.
\end{equation*}
If we set ${\mathbf f}_j = G({\mathbf s}_j) - {\mathbf s}_j =
{\mathbf s}_{j+1}- {\mathbf s}_j$ (where here ${\mathbf s}_{j+1}$ denotes the Picard
iteration)
we can compute a new vector ${\mathbf s}_{j+1}$ as
\begin{equation*}
{\mathbf s}_{j+1}=\mathbf{s}_{j}-G^{(j)}_{M,\lambda}\mathbf{f}_j.
\end{equation*}
Observe that $G^{(j)}_{M,\lambda}$ satisfies the following {\it
  multisecant} condition, see, e.g., \cite{fang2009two}, (when $\lambda=0$)
\begin{equation*}
	G^{(j)}_{M,\lambda}{\Delta}^2S_{j-m_j}^{(m_j)}=[\Delta {\mathbf s}_{j-m_j},\ldots, \Delta {\mathbf s}_{j-1}].
\end{equation*}
It is interesting to notice that, when $\lambda \neq 0$, we obtain a
class of {\it regularized projection methods}, that do not yet seem to have been  fully
investigated in the literature.

\vskip 2mm
For the sake of simplicity,  we did not present here the topological approaches of
Section \ref{sec:topological_approaches}, but the preceding algorithms can be easily
modified for these transformations.

\section{Anderson-type Mixing (ATM) methods} \label{sec:anderson_type}

\textit{Anderson Acceleration} ({AA}) (also known as \textit{Anderson Mixing}) is a technique originally presented in
\cite{ander} for solving systems of nonlinear equations written as
$F({\mathbf s})= G({\mathbf s}) - {\mathbf s}= 0$. In this section, we  generalize the basic
version of {AA} as given by Walker and Ni \cite{WalkerNi2011}  or by Higham and Strabi\'c \cite{a4}.  The
main idea of this generalization is that a procedure
similar to Anderson Acceleration can be built up with any of the
Shanks transformations.  We will name such methods {\it
Anderson-type Mixing} (ATM)  to emphasize the fact that, as it will be explained, these methods use a
Continuous-Updating scheme which {\it mixes} information coming out from
two different sequences.

Indeed, in the framework of the Continuous-Updating scheme presented
in Section \ref{sec:continuos_updating}, two different sequences are
generated, i.e., the continuously updated sequence $(\mathbf{s}_j)$ on
the one hand, and the sequence $(G(\mathbf{s}_j))$ on the other.
The main feature of
the Anderson-Mixing strategy is that it combines
the information coming from these two sequences in order to obtain a
better acceleration procedure. We will prove that it coincides
with a quasi-Newton strategy.
Since, in this case, the sequence $(\mathbf{s}_j)$ is not generated by a
\textit{fixed point iteration}, we also consider the
sequence $(\mathbf{f}_j)$, where $\mathbf{f}_j=
	G(\mathbf{s}_j)-\mathbf{s}_j= \mathbf{g}_j-\mathbf{s}_j$, that do not coincide with the
sequence $(\Delta \mathbf{s}_j)$.

 Algorithm \ref{ATMtable_1} shown below is a prototype version of {\it Anderson-type Mixing} method where
we define:
\begin{equation*}
  F_{j-m_j}^{(m_j)} \equiv [\mathbf{f}_{j-m_j},\ldots,  \mathbf{f}_{j-1}]
  \end{equation*}
  and use the previous notation
  $ S_{j-m_j}^{(m_j)} \equiv [\mathbf{s}_{j-m_j},\ldots,\mathbf{s}_{j-1}]$.

\vskip 2mm

\begin{algorithm}[H]
	\caption{The Anderson-Type Mixing (ATM) method.} \label{ATMtable_1}
	\DontPrintSemicolon
	\KwIn{Choose  $m \in \mathbb{N}, m \geq 1$, $\beta \in \mathbb{R}$, ${\mathbf s}_0 \in \mathbb{R}^p$.}
	Compute ${\mathbf f}_0=G({\mathbf s}_0)-{\mathbf s}_0$ and ${\mathbf s}_1={\mathbf s}_0+\beta {\mathbf f}_0$
 \;
	\For{$j=1,2,\ldots$}{
		{Compute $\mathbf{f}_j= G(\mathbf s_j)-\mathbf{s}_j$}\;
		{Set $m_j=\min(m,j)$}\;
		Set $ \Delta {S}_{j-m_j}^{(m_j)}=[ \Delta \mathbf{s}_{j-m_j}, \ldots, \Delta \mathbf{s}_{j-1}]$ and $ \Delta {F}_{j-m_j}^{(m_j)}=[ \Delta \mathbf{f}_{j-m_j}, \ldots, \Delta \mathbf{f}_{j-1}]$\;
		Compute $\boldsymbol{\theta}^{(j)} \in \mathbb{R}^{m_j}$\;
        Compute
		 $\overline{\mathbf s}_j= {\mathbf s}_j- \Delta {S}_{j-m_j}^{(m_j)}\boldsymbol{\theta}^{(j)}$   and
		 $\overline{\mathbf f}_j= {\mathbf f}_j- \Delta {F}_{j-m_j}^{(m_j)} \boldsymbol{\theta}^{(j)}$\;
		Set ${\mathbf s}_{j+1}=\overline{\mathbf s}_j+\beta \overline{\mathbf f}_j$  \; \label{algline:AAupdate}
		
	}
\end{algorithm}

\vskip 2mm
The scalar $\beta$,  usually a fixed positive value with $0 <\beta \leq 1$,   is called
\textit{mixing} or
\textit{damping} parameter. It is also possible to change it at each cycle, and it can be used to improve
convergence.
A common choice is to take $\beta=1$. In this case, since $\mathbf{g}_j= G(\mathbf{s}_j) =
\mathbf{s}_j + \mathbf{f}_j$ we can define
\begin{equation*}
  G_{j-m_j}^{(m_j)}=[\mathbf{g}_{j-m_j},\ldots,  \mathbf{g}_{j-1}]=
   {S}_{j-m_j}^{(m_j)} +   {F}_{j-m_j}^{(m_j)} .
  \end{equation*}
  By denoting $\overline{\mathbf{g}}_j=
\overline{\mathbf{s}}_j + \overline{\mathbf{f}}_j =
\mathbf{g}_{j} - \Delta {G}_{j-m_j}^{(m_j)}\boldsymbol{\theta}^{(j)}$, the new iterate can be simply computed as
${\mathbf s}_{j+1}=\overline{\mathbf g}_j$. This is the so-called
\textit{undamped iterate}.

Let us point out that  Line \ref{algline:AAupdate} in Algorithm \ref{ATMtable_1} can be alternatively written as
\begin{equation} \label{eq:general_theta_AA_UPDATE}
	\mathbf{s}_{j+1}={\mathbf{s}}_j- (-\beta \mathbf{f}_j+(\Delta {S}_{j-m_j}^{(m_j)}+\beta \Delta {F}_{j-m_j}^{(m_j)})\boldsymbol{\theta}^{(j)}),
\end{equation}
and that different choices of $\boldsymbol{\theta}^{(j)}$ give rise to
different ATMs. Some particular cases are described in the sequel.

The original {AA}
is obtained when
\begin{equation} \label{eq:rre_anderson}
	\boldsymbol{\theta}^{(j)}=\underset{\boldsymbol{\eta} \in \mathbb{R}^{m_j}}{\arg \min} \|\mathbf{f}_j - \Delta {F}_{j-m_j}^{(m_j)} \boldsymbol{\eta} \|^2,
\end{equation}
that is, assuming that the columns of $ \Delta {F}_{j-m_j}^{(m_j)} $ are linearly independent,
\begin{equation*}
	\boldsymbol{\theta}^{(j)}=
((\Delta {F}_{j-m_j}^{(m_j)})^T \Delta {F}_{j-m_j}^{(m_j)})^{-1}
(\Delta {F}_{j-m_j}^{(m_j)})^T\mathbf{f}_j .
\end{equation*}

\begin{rem}
It is interesting to observe that defining  $\boldsymbol{\theta}^{(j)}$'s as
	\begin{equation*}
	\boldsymbol{\theta}^{(j)}=\underset{\boldsymbol{\eta} \in \mathbb{R}^{m_j}}{\arg \min} \|\Delta \mathbf{s}_j -
	\Delta^2 S_{j-m_j}^{(m_j)}\boldsymbol{\eta}\|^2 \,
	\end{equation*}
	i.e., using \eqref{eq:RRE_strategy}, would be a good choice if the sequence $(\mathbf{s}_j)$ is
close to the Shanks kernel.  Instead, in the original AA the derivation of the
	$\boldsymbol{\theta}^{(j)}$  using \eqref{eq:rre_anderson} could be interpreted as an implicit assumption
that the sequence  $(\mathbf{f}_j)$ is \textit{closer} to the Shanks kernel than the sequence $(\mathbf{s}_j)$.
\end{rem}

  From \eqref{eq:general_theta_AA_UPDATE}, we have
\begin{equation}\label{eq:quasi_newton_variable_beta}
\mathbf{s}_{j+1}={\mathbf{s}}_j- (-\beta I +(\Delta S_{j-m_j}^{(m_j)}+\beta\Delta F_{j-m_j}^{(m_j)})((\Delta F_{j-m_j}^{(m_j)})^T\Delta F_{j-m_j}^{(m_j)})^{-1}(\Delta F_{j-m_j}^{(m_j)})^T )\mathbf{f}_j,
\end{equation}
as also observed in \cite{fang2009two,brzs}.

Formula \eqref{eq:general_theta_AA_UPDATE} highlights the connections between Anderson Mixing and quasi-Newton methods.
Indeed, in this case, defining

\begin{equation*}
	H_j^{(\beta)}= - \beta I +(\Delta S_{j-m_j}^{(m_j)}+ \beta \Delta F_{j-m_j}^{(m_j)})((\Delta F_{j-m_j}^{(m_j)})^T\Delta F_{j-m_j}^{(m_j)})^{-1}(\Delta F_{j-m_j}^{(m_j)})^T ,
\end{equation*}
we can write
 \begin{equation*}
 	\mathbf{s}_{j+1}={\mathbf{s}}_j- H_j^{(\beta)}\mathbf{f}_j,
 \end{equation*}
with $H_j^{(\beta)}$ satisfying the multisecant condition
$H_j^{(\beta)}\Delta F_{j-m_j}^{(m_j)}=\Delta S_{j-m_j}^{(m_j)}$. In
the next section we will fully make use of this idea: by introducing a
\textit{stabilization procedure} to overcome problems connected to the
ill-conditioning of the matrix $(\Delta F_{j-m_j}^{(m_j)})^T\Delta
F_{j-m_j}^{(m_j)}$, it is possible to prove the local linear convergence of the AA method.

As indicated in the previous Section, it is also possible to define a
  {\it Periodic Anderson-Type Mixing} method whereby acceleration steps are interspersed into
  linear updates  at regular intervals.
  Fixing the period
$\mu \in \mathbb{N}$, with $\mu \geq 1$, an
Anderson-type update is made each $\mu$ iterations.
In between these updates, when $\mu > 1$,
the iterates are computed simply as a  linear mixing ${\mathbf s}_{j+1}={\mathbf s}_{j}+\beta {\mathbf f}_{j}$,
where $\beta$ is the mixing parameter ($\beta=1$ corresponds to Picard's iteration).
Clearly, when  $\mu = 1$  Algorithm \ref{ATMtable_P}  coincides with
Algorithm \ref{ATMtable_1}.

\vskip 2mm

\begin{algorithm}[H]
	\caption{The Periodic Anderson-Type Mixing method.} \label{ATMtable_P}
	\DontPrintSemicolon
	\KwIn{Choose  $m, \mu \in \mathbb{N}, \, m, \mu \geq 1$, $\beta \in \mathbb{R}$, ${\mathbf s}_0 \in \mathbb{R}^p$.}
	Compute ${\mathbf f}_0=G({\mathbf s}_0)-{\mathbf s}_0$ and ${\mathbf s}_1={\mathbf s}_0+\beta {\mathbf f}_0$
 \;
	\For{$j=1,2,\ldots$}{
		Compute $\mathbf{f}_j= G(\mathbf s_j)-\mathbf{s}_j$\;
\eIf{$(j+1) \!\!\mod \mu = 0$}{
        Compute ${\mathbf s}_{j+1}$ using steps 4 to 8 of
        Algorithm \ref{ATMtable_1}  {\it (Anderson-type update)} \;}
{Compute ${\mathbf s}_{j+1}={\mathbf s}_j+\beta {\mathbf f}_j$
{\it (linear mixing update)}
}
	}
\end{algorithm}

\vskip 2mm
  It is important to underline that the values of $\mu$ and $m$ can be chosen independently.
  However, when $\mu \ge 3$ and  we choose $m=\mu-2$, then
 the terms used for computing the Anderson-type update are only those  terms obtained by the linear mixing update, and therefore in this situation
  Algorithm \ref{ATMtable_P} proceeds as a RM method of
Algorithm \ref{RMtable}, with a different restarting formula.
It must also be noticed that Algorithm \ref{ATMtable_P}  with $\boldsymbol{\theta}^{(j)}$ computed as in
 \eqref{eq:rre_anderson},  is exactly the Periodic Pulay method \cite{A} (compare also with \eqref{eq:quasi_newton_variable_beta}).
Interleaving Anderson Acceleration with
fixed point iterations for improving the global convergence properties, but not
necessarily the speed, has been recognized before in the physics literature
as can be seen from the related discussion and the references in \cite{A}.
This idea is somewhat similar also to the A2DR
(Anderson accelerated Douglas--Rachford) algorithm proposed in  \cite{B}.

\vskip 0.2cm
To start the derivation of the new ATMs, we observe that a  possible generalization for the derivation of the $\boldsymbol{\theta}^{(j)}$
 can be obtained by using the coupled sequences defined in Section \ref{coupled}, that is by taking
\begin{equation}\label{eq:coupledtheta}
	\boldsymbol{\theta}^{(j)}=
(Y^T \Delta {C}_{j-m_j}^{(m_j)})^{-1}
Y^T\mathbf{c}_j .
\end{equation}
If we take $\mathbf{c}_j= \mathbf{f}_j$, for all $j$, and $Y = \Delta {C}_{j-m_j}^{(m_j)}=
\Delta {F}_{j-m_j}^{(m_j)}$ we recover the AA choice for  $\boldsymbol{\theta}^{(j)}$.
It is easy to see that taking into account the transformations defined at the beginning of Section
\ref{sec:extr_techniques},
if we consider the extrapolated vector
$
\widetilde{\boldsymbol{t}}_{j-m_j}^{(m_j+1)}=
{\mathbf s}_j- \Delta {S}_{j-m_j}^{(m_j)}\boldsymbol{\theta}^{(j)}$ we recover
exactly  the $\overline{\mathbf{s}}_j$'s computed
in Algorithms \ref{ATMtable_1} and
\ref{ATMtable_P}. If we consider the same $\boldsymbol{\theta}^{(j)}$, in the same
relation, and by using, as sequence to be extrapolated the coupled one $(\mathbf{f}_j)$,
we obtain $\overline{\mathbf{f}}_j$.

Another additional generalization can be made by considering, as in problem \eqref{eq:RRE_problem_prec}
 of Section \ref{sec:extr_techniques}, a different metric in the evaluation of the norm, and also a regularization parameter
$\lambda$. We consider the problem
\begin{equation}\label{eq:AA_type_coefficients_problem}
\boldsymbol{\theta}^{(j)}_{M,\lambda}=\underset{\boldsymbol{\eta} \in \mathbb{R}^{m_j}}{\arg \min}
\left(\|\mathbf{c}_j - \Delta {C}_{j-m_j}^{(m_j)}\boldsymbol{\eta}\|_M^2+\lambda \|\boldsymbol{\eta}\|^2\right).
\end{equation}
The solution  is
\begin{equation}\label{eq:dopo34}
\boldsymbol{\theta}^{(j)}_{M,\lambda}=((\Delta C_{j-m_j}^{(m_j)})^TM\Delta C_{j-m_j}^{(m_j)}+\lambda I)^{-1}(\Delta C_{j-m_j}^{(m_j)})^TM
\mathbf{c}_j.
\end{equation}

By taking in \eqref{eq:dopo34} $\mathbf{c}_j= \mathbf{f}_j$ and $M=I$,  that is by introducing only a $\ell_2$-regularization term to the original AA problem, we obtain a method that we call   \textit{Regularized
Anderson Acceleration} (in short RAA).

If we take $M=YY^T$ and $\lambda=0$, it is possible to see that $\boldsymbol{\theta}^{(j)}_{M,\lambda}$ in \eqref{eq:dopo34} can be obtained, alternatively, as the solution of the linear
 system
\begin{equation*}
(\Delta C_{j-m_j}^{(m_j)})^TY(Y^T\Delta C_{j-m_j}^{(m_j)}\boldsymbol{\theta}^{(j)}_{YY^T \!,0}- Y^T \mathbf{c}_j)=0,
\end{equation*}
which correspond exactly to \eqref{eq:coupledtheta}, assuming that
$\mbox{\rm rank}( (\Delta C_{j-m_j}^{(m_j)})^TY )= m_j$.

The ATMs methods can thus be obtained by considering the coupled
sequence $(\mathbf{c}_j)= (\mathbf{f}_j)$ fixed, and changing the
matrix $Y$. The following particular cases are of interest:

\begin{enumerate}

\item ATM-RRE: $Y=[\Delta^2 {\mathbf s}_{j-m_j},\ldots, \Delta^2 {\mathbf s}_{j-1}]= \Delta^2 S_{j-m_j}^{(m_j)}\in {\mathbb R}^{p \times m_j}$ corresponds to a method in the style of the RRE. For this choice, since we also need the knowledge  of the vector
  ${\mathbf s}_{j+1}$ we have to edit slightly Algorithm \ref{ATMtable_1}
    by beginning the loop (line 2) with $j=2$ and by adding before it
    the  computation of ${\mathbf s}_2={\mathbf s}_1+\beta {\mathbf f}_1$.
    Modifications that take this into account must also be made in Algorithm \ref{ATMtable_P}.
    The choice $Y=[\Delta^2 {\mathbf f}_{j-m_j},\ldots, \Delta^2 {\mathbf f}_{j-1}]= \Delta^2 F_{j-m_j}^{(m_j)}\in {\mathbb R}^{p \times m_j}$ is also possible.

\item ATM-MPE: $Y=[\Delta {\mathbf s}_{j-m_j},\ldots, \Delta {\mathbf s}_{j-1}]
= \Delta S_{j-m_j}^{(m_j)}\in {\mathbb R}^{p \times m_j}$ or $Y=[{\mathbf f}_{j-m_j},\ldots, {\mathbf f}_{j-1}]=F_{j-m_j}^{(m_j)}\in {\mathbb R}^{p \times m_j}$  leads to two methods in the style of the MPE;

\item ATM-MMPE: $Y=[{\mathbf y}_1,\ldots,{\mathbf y}_{m_j}]\in {\mathbb R}^{p \times m_j}$ which leads to an ATM in the style of the MMPE.

\item ATM-TEA: suitably modifying the structure of
Algorithms \ref{ATMtable_1} and \ref{ATMtable_P}, it is possible to use a topological approach (see
  Section \ref{sec:topological_approaches}) to obtain the coefficients
  $\boldsymbol{\theta}^{(j)}_{M,\lambda}$. As in Section
  \ref{sec:continuos_updating} we omit all the details for the sake of
  brevity.
\end{enumerate}

Before concluding this section, we point out that  the introduction of an $\ell_2$-regularization term for AA has
already been studied in the recent papers \cite{B,C,ander1}, and that \eqref{eq:AA_type_coefficients_problem} represents
a generalization to the ATM methods of the $\ell_2$-regularization approach for AA. In Section \ref{sec:numerical_res},
for the particular AA case,  we will propose and experimentally analyze the choice of the regularization parameter
$\lambda$ using the \textit{Generalized Cross Validation} \cite{GCVgolub}. This choice represents a major difference
with the above mentioned works, where the choice of the regularization parameter is made adaptively based on quantities
related to the most recent iterates (see, for example, \cite[eq. (3.4)]{B} and \cite[eq. (3)]{C}). Sections
\ref{sec:Stabilized_AA} and \ref{sec:stabilization} below
further justify/clarify the introduction of an $\ell_2$-regularization strategy.

\subsection{Stabilized AA} \label{sec:Stabilized_AA}

The aim of this Section is to present an algorithm  which can be {viewed} as a \textit{stabilized}
version of the {AA} method. In particular, in this new version of {AA},
  a check on the linear independence of the vectors $\Delta
  \mathbf{f}_d$ is performed (Lines \ref{alg:lineli}
  -\ref{alg:linelif}): the residual difference $\Delta \mathbf{f}_{d}$
  is discarded if its projection $\widehat{\mathbf{f}}_{d}$
  onto the orthogonal of the previously computed residual differences
  is \textit{close} to the null vector, i.e., if it results in a vector
  of sufficiently small norm when compared to the original one (see
  Section \ref{sec:stabilization} for further details). It is interesting to note that when, in Algorithm \ref{alg:ATM_2},
  we choose $m=1$ (and likely for small values of $m$) the introduced stabilization procedure is not required
  and Algorithm \ref{alg:ATM_2} coincides with the classic AA scheme (compare, in this case,
  \eqref{eq:quasi_newton_variable_beta} and the update at Line \ref{algline:update} in Algorithm~\ref{alg:ATM_2}).

\vskip 0.2cm
\begin{algorithm}[H]
\caption{Stabilized Anderson Acceleration.} \label{alg:ATM_2}
	\DontPrintSemicolon
	\KwIn{Choose  $m \in \mathbb{N}, m \geq 1$, $\beta \in \mathbb{R}$,  ${\mathbf s}_0 \in \mathbb{R}^p$ and $\tau>1$.}
	Compute ${\mathbf f}_0=G({\mathbf s}_0)-{\mathbf s}_0$ and ${\mathbf s}_1={\mathbf s}_0+\beta {\mathbf f}_0$
 \;
	\For{$j=1,\ldots$}{
		Set $m_j=\min(m,j).$ \;
		Compute $\mathbf {f}_j= G(\mathbf{s}_j)-\mathbf{s}_j$  \;
		Compute $\widehat {\mathbf{f}}_{j-m_j}= {\Delta \mathbf{f}_{j-m_j}}$ \;
Set  $P_{{j-m_j}}= ({\widehat {\mathbf{f}}_{{j-m_j}}\widehat {\mathbf{f}}_{{j-m_j}}^T})/({\widehat {\mathbf{f}}_{{j-m_j}}^T\widehat {\mathbf{f}}_{{j-m_j}}})$\;
		\For{$d=j-m_j+1,\ldots, j-1$  \label{alg:lineli}}{
		Set $Q_{{j-m_j}}^{{d-1}}= \sum_{i=j-m_j}^{d-1} P_{i}$ \;
		Compute $\widehat {\mathbf{f}}_{d}=(I-Q_{{j-m_j}} ^{{d-1}}){\Delta \mathbf{f}_{d}}$ \;
		\eIf{$\|\widehat {\mathbf{f}}_{d}\|\tau \geq \|{\Delta \mathbf{f}_{d}}\|$}{
		Set $P_{d}= ({\widehat {\mathbf{f}}_{d}\,\widehat {\mathbf{f}}_{d}^T})/({\widehat {\mathbf{f}}_{d}^T \,\widehat {\mathbf{f}}_{d}})$\;
	}{
		Set $\widehat {\mathbf{f}}_{d}= \boldsymbol{0}$ \;
		Set $P_{d}=\boldsymbol{0}$\;}
	} \label{alg:linelif}
	 Set $\mathcal{I}_j=   \{k_1, \ldots, k_{\widehat{m}_j}\} \subseteq \{ {j-m_j}, \ldots, {j-1} \}$ the set of indices such that $\widehat{\mathbf{f}}_{k_1}, \cdots, \widehat{\mathbf{f}}_{k_{\widehat{m}_j}}$ are non null vectors\; \label{algline:definitionhatfkd}
		Set $\Delta F_{\mathcal{I}_j}= [\Delta \mathbf{f}_{k_1},\ldots,\Delta \mathbf{f}_{k_{\widehat{m}_j}} ], \; \Delta S_{\mathcal{I}_j}= [\Delta \mathbf{s}_{k_1},\ldots,\Delta \mathbf{s}_{k_{\widehat{m}_j}} ]$ \;
		Set $H_j^{(\beta)}=  \big[- \beta I +(\Delta S_{\mathcal{I}_j}+ \beta \Delta F_{\mathcal{I}_j})((\Delta F_{\mathcal{I}_j})^T\Delta F_{\mathcal{I}_j})^{-1}(\Delta F_{\mathcal{I}_j})^T \big  ]$\; \label{alg:qnMtrices}
		Compute ${\mathbf s}_{j+1}={\mathbf s}_j-H_j^{(\beta)}{\mathbf f}_j$  \;\label{algline:update}
		}
\end{algorithm}

\subsubsection{Local convergence}

There already exist in the literature  different proofs of the  local convergence for  the stabilized
versions of AA, see for example \cite{RthesisQC,GSsnep,B,C,RSDIIS}.
In principle, our convergence analysis can be obtained using
ideas and techniques from \cite[Sec. 4.2]{RthesisQC}, but we prefer to present here a full detailed proof. The reasons to
present such a detailed proof can be mainly summarized as follows: a) our derivation is not completely analogous
to that in \cite{RthesisQC}: simplifying some arguments, we are able to obtain slightly more general results than those
presented in \cite[Sec. 4.2]{RthesisQC} (the interested reader can compare our Theorem~\ref{th:AAconvergence_theorem}
with \cite[Th. 4.10]{RthesisQC}) ; b) our analysis does not require the contractivity or non-expansivity of the fixed point
map $G$, a major difference if compared to what has been proved in \cite{C,B}; c) our proof of convergence holds for
every mixing parameter $\beta \in \mathbb{R}$ shedding further light on the significance and the relevance of the
parameter $\beta$ in the {AA} procedure: it can be interpreted as a scaling factor of the initial Jacobian approximation
(see Theorem~\ref{th:AAconvergence_theorem}); d) when $m=1$, since Algorithm~\ref{alg:ATM_2} coincides with
the classic AA scheme (see the beginning of Section~\ref{sec:Stabilized_AA}), we obtain, as a by-product of our analysis,
an alternative proof of that given in \cite[Sec. 2.3]{toth2015convergence} for the convergence of the classic AA with
$m=1$ without assuming, once more, any contractivity of the fixed point map $G$.
We consider the function
$F(\mathbf{s})=G(\mathbf{s})-\mathbf{s}$, and  we
made the following assumption:

\begin{assumptions} \label{assumption:diff}
$F: \mathbf{R}^{n} \rightarrow \mathbf{R}^{n} $ is
  differentiable in a open convex set $E\subseteq \mathbb{R}^n$ and
  there exists $\mathbf{s}^{*} \in E$ such that
  $\mathbf{f}^{*}= F(\mathbf{s}^*)=\boldsymbol{0}$.  Moreover,
   $J=F'(\mathbf{s}^*)$  is invertible and   for all $\mathbf{s} \in  E$
  we have
  \begin{equation*}
    \|F'(\mathbf{s})-F'(\mathbf{s}^*)\| \leq L \|\mathbf{s}-\mathbf{s}^*\|.
  \end{equation*}
 The above assumption implies that,
  \begin{equation*}
    \|F(\mathbf{u})-F(\mathbf{v})-J(\mathbf{u}-\mathbf{v})\|\leq L \|\mathbf{u}-\mathbf{v}\|\max\{\|\mathbf{u}-\mathbf{s}^*\|,\|\mathbf{v}-\mathbf{s}^*\|\},
  \end{equation*}
  for all $\mathbf{u},\mathbf{v} \in E$ and that there exists $U_{\kappa}(\mathbf{s}^*):=\{\mathbf{u} \in \mathbb{R}^n \;:\;  \|\mathbf{u}- \mathbf{s}^{*}\| \leq \kappa\}$ $\hbox{ s.t., for some } \rho>0$,
	\begin{equation*}
	\frac{1}{\rho}\|\mathbf{u}-\mathbf{v}\| \leq \|F(\mathbf{u})-F(\mathbf{v})\| \leq \rho \|\mathbf{u}-\mathbf{v}\|.
	\end{equation*}
\end{assumptions}

In the remainder of this section we use the notations introduced in Algorithm \ref{alg:ATM_2}.

\begin{lemma} \label{lem:Multisecant_stabilized}
	The matrices $H_j^{(\beta)}$ {(defined at Line \ref{alg:qnMtrices} of Algorithm \ref{alg:ATM_2})} satisfy the multisecant condition
	\begin{equation*}
		H_j^{(\beta)}\Delta F_{\mathcal{I}_j}=\Delta S_{\mathcal{I}_j}
	\end{equation*}
\end{lemma}
\noindent{\bf Proof:}
The  proof is by direct verification.
\fin
\begin{lemma} \label{lemma:inductive_multisecan}
	$H_j^{(\beta)}$ can be computed  {recursively} from $H_j^{0}=-\beta I$ using

	\begin{equation*}
		H_j^{d}=H_{j}^{d-1}+\frac{(\Delta \mathbf{s}_{k_d}-H_j^{d-1}\Delta \mathbf{f}_{k_d})\widehat {\mathbf{f}}_{k_d}^T}
{\widehat{\mathbf{f}}_{k_d}^T{\Delta \mathbf{f}_{k_d}}} \hbox{ for } d=1,\ldots, \widehat{m}_j
	\end{equation*}
	with $H_j^{\widehat{m}_j}=H_j^{(\beta)}$  {(see Line \ref{algline:definitionhatfkd} in Algorithm \ref{alg:ATM_2} for the definitions of $\widehat{\mathbf{f}}_{k_d}$)}. In particular, for all $d=1,\ldots, \widehat{m}_j$,
        we have: $H_j^{d}\Delta\mathbf{f}_{k_p}=\Delta\mathbf{s}_{k_p}$ for $p=1,\ldots,d$.
\end{lemma}

\noindent{\bf Proof:} Define $Z \in \mathbb{R}^{n \times n-\widehat{m}_j}$ as
a basis for ${\rm{span}}(\Delta F_{\mathcal{I}_j})^{\perp} $. From the definition
of $H_j^{(\beta)}$ we have $H_j^{(\beta)}Z=-\beta Z$ and
$H_j^{(\beta)}\Delta F_{\mathcal{I}_j}=\Delta S_{\mathcal{I}_j}$. To
prove the theorem, we will prove (by induction) that
$H_{j}^{\widehat{m}_j}$ satisfies the same relations. For $d=1$ we
have $\displaystyle H_j^1=H_j^0+\frac{(\Delta \mathbf{s}_{k_1}-H_j^{0}\Delta
  \mathbf{f}_{k_1})\widehat {\mathbf{f}}_{k_1}^T}{\widehat {\mathbf{f}}_{k_1}^T{\Delta \mathbf{f}_{k_1}}}$ and hence
$H_j^1\Delta \mathbf{f}_{k_1}=\Delta \mathbf{s}_{k_1}$. Suppose now
the assumption true for $d=\ell$. By definition we have that
$H_j^{\ell+1}\Delta \mathbf{f}_{k_{\ell+1}}=\Delta
\mathbf{s}_{k_{\ell+1}}$ and $H_j^{\ell+1}\Delta
\mathbf{f}_{k_{p}}=\Delta \mathbf{s}_{k_{p}}$ for all $p=1,\ldots,\ell$
since $\widehat {\mathbf{f}}_{k_{\ell+1}}\perp \Delta\mathbf{f}_{k_{p}}$. Finally, since $${\rm{span}}(\widehat {\mathbf{f}}_{k_{1}},\ldots, \widehat {\mathbf{f}}_{k_{\widehat{m}_j}})=
{\rm{span}}({\Delta\mathbf{f}_{k_{1}}},\ldots,
{\Delta\mathbf{f}_{k_{\widehat{m}_j}}} ), $$  implies that $Z$
is also a basis for $ {\rm{span}}(\widehat {\mathbf{f}}_{k_{1}},\ldots, \widehat {\mathbf{f}}_{k_{\widehat{m}_j}})^{\perp}$, we have $H_j^{\widehat{m}_j}Z=-\beta Z$. {The result} follows observing that,
since $[\Delta F_{\mathcal{I}_j},Z]$ is invertible, the equation
$B[\Delta F_{\mathcal{I}_j},Z]=[\Delta S_{\mathcal{I}_j},-\beta Z]$
has a unique solution.  \fin

Observe that, as already pointed out in \cite{zhangboyd}, Lemma \ref{lemma:inductive_multisecan}
highlights the connections between the Jacobian approximations produced by the \textit{Bad (or type-II) Broyden} update \cite{badBroyden}
and the matrices produced by AA.

\begin{lemma} \label{lemma:recursiveH_overline}
Let us define $\widehat {\mathbf{s}}_{k_1}= \Delta \mathbf{s}_{k_1}$ and for $d=2,\ldots,\widehat{m}_j$ define  $\widehat { \mathbf{s}}_{k_d}= \Delta \mathbf{s}_{k_d}-H_j^{d-1}Q_{k_1}^{k_{d-1}} \Delta \mathbf{f}_{k_d}$ being $Q_{k_1}^{k_{d-1}}=\sum_{p=1}^{d-1}({\widehat {\mathbf{f}}_{k_p}\widehat {\mathbf{f}}_{k_p}^T}/{\widehat {\mathbf{f}}_{k_p}^T\widehat {\mathbf{f}}_{k_p}})$.
Then $H_j^{(\beta)}$ can be computed  {recursively} from $H_j^{0}=-\beta I$ using
	\begin{equation*}
	H_j^{d}=H_{j}^{d-1}+\frac{(\widehat {\mathbf{s}}_{k_d}-H_j^{d-1}\widehat {\mathbf{f}}_{k_d})\widehat {\mathbf{f}}_{k_d}^T}{\widehat {\mathbf{f}}_{k_d}^T\widehat {\mathbf{f}}_{k_d}} \hbox{ for } d=1,\ldots, \widehat{m}_j
	\end{equation*}
	with $H_j^{\widehat{m}_j}=H_j^{(\beta)}$. In particular, for
        all $d=1,\ldots, \widehat{m}_j$, we have:
        $H_j^{d}\widehat {\mathbf{f}}_{k_p}=\widehat {\mathbf{s}}_{k_p}$
        for $p=1,\ldots,d$.
\end{lemma}
\noindent{\bf Proof:} The proof follows  {from}  the definition of $H_j^{d}$,
and observing
that $$\widehat {\mathbf{f}}_{k_d}=(I-Q_{k_1}^{k_{d-1}}){\Delta\mathbf{f}_{k_d}}
\Rightarrow
\widehat {\mathbf{f}}_{k_d}^T\widehat {\mathbf{f}}_{k_d}= {\widehat {\mathbf{f}}_{k_d}}^T{\Delta\mathbf{f}_{k_d}}$$
(since $(I-Q_{k_1}^{k_{d-1}})$ is a projector) and that $
\widehat {\mathbf{s}}_{k_d}
-H_j^{d-1}\widehat {\mathbf{f}}_{k_d}={\Delta\mathbf{s}_{k_d}}
-H_j^{d-1}{\Delta\mathbf{f}_{k_d}}.$\fin
\begin{lemma}\label{lem:bounded_deterioration}
  Suppose that $\mathbf{s}_{k_d}, \mathbf{s}_{k_d+1} \in U_{\kappa}(\mathbf{s}^{*})$ for all $d=1,\ldots,\widehat{m}_j$.
  Then, the following inequality is satisfied
	\begin{equation*}
		\|\widehat {\mathbf{s}}_{k_d}-J^{-1}\widehat {\mathbf{f}}_{k_d}\|\leq C \|{\Delta \mathbf{f}_{k_d}} \|
\sum_{p=1}^{d}n_{k_p}^{k_p+1} (2 \tau)^{p-d},
	\end{equation*}
	where $C= \|J^{-1}\| L \rho $ and $n_{k_p}^{k_p+1}= \max \{\|\mathbf{s}_{k_p+1}-\mathbf{s}^*\|,\|\mathbf{s}_{k_p}-\mathbf{s}^*\| \}$.
\end{lemma}
\noindent{\bf Proof:}
	For $d=1$ we have
	$$\|\widehat {\mathbf{s}}_{k_1}-J^{-1}\widehat {\mathbf{f}}_{k_1}\|=\|{\Delta \mathbf{s}_{k_1}}-J^{-1}
{\Delta \mathbf{f}_{k_1}}\|\leq C \|\Delta \mathbf{f}_{k_1}\|n_{k_1}^{k_1+1},$$
where  the last inequality follows  {from} Assumption \ref{assumption:diff}.
Suppose now the assumption true for $d=\ell$. To prove the statement for $d=\ell+1$ we have
	\begin{equation*}
	\begin{split}
	& \|\widehat {\mathbf{s}}_{k_{\ell+1}}-J^{-1}\widehat {\mathbf{f}}_{k_{\ell+1}}\| \leq \|\Delta \mathbf{s}_{k_{\ell+1}}-J^{-1}\Delta \mathbf{f}_{k_{\ell+1}}\|+\|H_j^{\ell}Q_{k_1}^{k_\ell}\Delta \mathbf{f}_{k_{\ell+1}}-J^{-1}Q_{k_1}^{k_\ell}\Delta \mathbf{f}_{k_{\ell+1}}  \| \\\
	& \leq C\|\Delta \mathbf{f}_{k_{\ell+1}}\| n_{k_{\ell+1}}^{k_{\ell+1}+1}+\sum_{p=1}^{\ell}\frac{\|H_j^{\ell}\widehat {\mathbf{f}}_{k_p}-J^{-1}\widehat {\mathbf{f}}_{k_p}  \|}{\|\widehat {\mathbf{f}}_{k_p}\|}\|\Delta \mathbf{f}_{k_{\ell+1}}\| \\
	& =C\|\Delta \mathbf{f}_{k_{\ell+1}}\| n_{k_{\ell+1}}^{k_{\ell+1}+1}+\sum_{p=1}^{\ell}\frac{\|\widehat {\mathbf{s}}_{k_p}-J^{-1}\widehat {\mathbf{f}}_{k_p}  \|}{\|\widehat {\mathbf{f}}_{k_p}\|}\|\Delta \mathbf{f}_{k_{\ell+1}}\|  \\
	& \leq C\|\Delta \mathbf{f}_{k_{\ell+1}}\|(n_{k_{\ell+1}}^{k_{\ell+1}+1}+\tau \sum_{p=1}^{\ell}\sum_{h=1}^{p}n_{k_h}^{k_{h}+1}(2\tau)^{p-h})= C\|\Delta \mathbf{f}_{k_{\ell+1}}\|(n_{k_{\ell+1}}^{k_{\ell+1}+1}+\tau \sum_{p=1}^{\ell}n_{k_p}^{k_{p}+1}
h=0^{\ell-p}(2\tau)^{h})
	\end{split}
	\end{equation*}
	where, in the first inequality, we use the definition of $\widehat {\mathbf{s}}_{k_{\ell+1}}$, in the second inequality, we use
the definition of $Q_{k_1}^{k_\ell}$, in the first equality, we use the fact that
$H_j^{\ell}\,\widehat {\mathbf{f}}_{k_p}=\widehat {\mathbf{s}}_{k_p}$ for $p=1,\ldots,\ell$
(see Lemma \ref{lemma:recursiveH_overline}), and, in the last inequality, our  {induction} hypothesis. Finally, since
	\begin{equation*}
		\sum_{h=0}^{\ell-p}(2\tau)^h \leq \tau^{\ell-p}\sum_{h=0}^{\ell-p}2^h=\tau^{\ell-p}(2^{\ell-p+1}-1)\leq \tau^{\ell-p}2^{\ell-p+1},
	\end{equation*}
	we have that
	\begin{equation*}
		C\|\Delta \mathbf{f}_{k_{\ell+1}}\|(n_{k_{\ell+1}}^{k_{\ell+1}+1}+\tau \sum_{p=1}^{\ell}n_{k_p}^{k_{p}+1}\sum_{h=0}^{\ell-p}(2\tau)^{h})\leq C\|\Delta \mathbf{f}_{k_{\ell+1}}\| \sum_{p=1}^{\ell+1}n_{k_p}^{k_{p}+1}(2\tau)^{\ell+1-p}
	\end{equation*}
	which concludes the proof.\fin

\begin{lemma} \label{lemma:boundend_deterioration_complete}
	The following equality is satisfied
	\begin{equation*}
		H_j^{(\beta)}-J^{-1}=(-\beta I-J^{-1})(I-Q_{k_{1}}^{k_{\widehat{m}_j}})+\sum_{d=1}^{\widehat{m}_j}\frac{(\widehat {\mathbf{s}}_{k_d}-J^{-1}\widehat {\mathbf{f}}_{k_d})\widehat {\mathbf{f}}_{k_d}^T }{\widehat {\mathbf{f}}_{k_d}^T\widehat {\mathbf{f}}_{k_d}}.
	\end{equation*}
	Moreover, if $\mathbf{s}_{k_d}, \mathbf{s}_{k_d+1} \in U_{\kappa}(\mathbf{s}^{*})$ and $n_{k_d}^{k_d+1}\leq \varepsilon$ for all $d=1,\ldots, \overline{m}_{j}$, there exists a constant
$\alpha=\alpha(\tau,m,C)$ such that
	\begin{equation*}
		\sum_{d=1}^{\widehat{m}_j}\frac{\|(\widehat {\mathbf{s}}_{k_d}-J^{-1}\widehat {\mathbf{f}}_{k_d})\widehat {\mathbf{f}}_{k_d}^T\| }{\widehat {\mathbf{f}}_{k_d}^T\widehat {\mathbf{f}}_{k_d}}\leq \alpha \varepsilon.
	\end{equation*}
\end{lemma}
\noindent{\bf Proof:}
	The first part of the statement follows from direct computation using the fact that the vectors $\widehat {\mathbf{f}}_{k_d}$
are orthogonal (see also \cite[Lemma 4.17]{RthesisQC}). For the second part, observe that
	\begin{equation*}
	\begin{split}
	& \sum_{d=1}^{\widehat{m}_j}\frac{\|(\widehat {\mathbf{s}}_{k_d}-J^{-1}\widehat {\mathbf{f}}_{k_d})\widehat {\mathbf{f}}_{k_d}^T\| }{\widehat {\mathbf{f}}_{k_d}^T\widehat {\mathbf{f}}_{k_d}}\leq
	 \sum_{d=1}^{\widehat{m}_j}\frac{\|(\widehat {\mathbf{s}}_{k_d}-J^{-1}\widehat {\mathbf{f}}_{k_d})\| }{\|\widehat {\mathbf{f}}_{k_d}\|}  \\
	 & \leq C\sum_{d=1}^{\widehat{m}_j}\frac{\|{\Delta \mathbf{f}_{k_d}}\| }{\|\widehat {\mathbf{f}}_{k_d}\|}\sum_{p=1}^{d}n_{k_d}^{n_d+1} (2\tau)^{d-p}\leq \varepsilon C\tau \sum_{d=1}^{\widehat{m}_j}\sum_{p=1}^{d} (2\tau)^{d-p}  \\
	 & \leq \varepsilon C\tau \sum_{d=1}^{{m}}\sum_{h=0}^{d-1} (2\tau)^{h}\leq \varepsilon C m(2\tau)^{m},
	\end{split}
	\end{equation*}
	where, in the second inequality, we use Lemma \ref{lem:bounded_deterioration}, and, in the fourth one,
the fact that $\widehat{m}_j \leq m_j \leq m$ for all $j$.
 \fin

\begin{theorem}	\label{th:AAconvergence_theorem}
	Let $\mathbf{s}_0, \mathbf{s}_1, \ldots, $ be the iterates produced by Algorithm~\ref{alg:ATM_2} (Stabilized Anderson Acceleration).
Then, for all $q \in (0,1)$, there exists $\delta = \delta (q, \alpha)$, $\varepsilon(q, \alpha)$ such that if
	\begin{equation*}
		\|-\beta I-J^{-1}\| \leq \delta \hbox{ and } \|\mathbf{s}_0-\mathbf{s}^{*}\| \leq \varepsilon,
	\end{equation*}
	we have
	\begin{equation*}
		\mathbf{s}_{j+1} \in E \hbox{ and }  \|\mathbf{s}_{j+1}-\mathbf{s}^*\|\leq q \|\mathbf{s}_j- \mathbf{s}^*\|
	\end{equation*}
	for all $j \in \mathbf{N}$.
\end{theorem}
\noindent{\bf Proof:}
	For a fixed $q$, choose $\delta$ and $\varepsilon$ such that
	\begin{equation*}
		\|J^{-1}\|L\varepsilon + \rho (\delta + \alpha \varepsilon) < q
	\end{equation*}
 in a way that $U_\varepsilon(\mathbf{s}^*)\subseteq U_\kappa (\mathbf{s}^*) \subseteq E$ (where $\kappa$
 and $\alpha$ are the same as in  Lemma \ref{lemma:boundend_deterioration_complete}, and $\rho$ is the same as in
 Assumption \ref{assumption:diff}).
	For $j=0$ we have
	\begin{equation*}
		\begin{split}
		& \|\mathbf{s}_1-\mathbf{s}^{*}\| \leq \|\mathbf{s}_0+\beta\mathbf{f}_0- \mathbf{s}^*\|\leq \|\mathbf{s}_0 - \mathbf{s}^*-J^{-1}(\mathbf{f}_0-\mathbf{f}^*)\|+\|(-\beta I-J^{-1})(\mathbf{f}_0-\mathbf{f}^*)\|  \\
		& \leq \|J^{-1}\|\|J(\mathbf{s}_0-\mathbf{s}^*)-(\mathbf{f}_0-\mathbf{f}^*)\|+ \delta \|\mathbf{f}_0-\mathbf{f}^*\| \leq ( \|J^{-1}\|L \varepsilon+ \delta \rho ) \|\mathbf{s}_0 -\mathbf{s}^*\| \leq q \|\mathbf{s}_0 -\mathbf{s}^*\|\leq \varepsilon,
		\end{split}
	\end{equation*}
	which proves that $\mathbf{s}_1 \in U_{\varepsilon}(\mathbf{s}^{*})$. Assume now that, for all $j \geq 0$,
$\|\mathbf{s}_j-\mathbf{s}^*\|\leq q^{j} \|\mathbf{s}_0-\mathbf{s}^*\|$
and hence that $\mathbf{s}_j \in U_\varepsilon(\mathbf{s}^*)$. We have
	\begin{equation*}
		\begin{split}
		& \|\mathbf{s}_{j+1}-\mathbf{s}^{*}\|=\|\mathbf{s}_j-H_j^{(\beta)}\mathbf{f}_j-\mathbf{s}^{*}\| \\
		& \leq \|J^{-1}\|\|J(\mathbf{s}_j-\mathbf{s}^{*})-(\mathbf{f}_j-\mathbf{f}^{*})\|+\|H_j^{(\beta)}-J^{-1}\|\|\mathbf{f}_j-\mathbf{f}^{*}\| \\
		& \leq \|J\|^{-1}L\|\mathbf{s}_j-\mathbf{s}^*\|^2+\rho \|H_j^{(\beta)}-J^{-1}\|\|\mathbf{s}_j-\mathbf{s}^*\|   \\
		& \leq (\|J\|^{-1}L q^{j} \varepsilon+ \rho (\delta +\alpha \varepsilon))\|\mathbf{s}_j-\mathbf{s}^*\| \leq q \|\mathbf{s}_j-\mathbf{s}^*\| ,
		\end{split}
	\end{equation*}
	where, in the last inequality, we use our {induction} hypothesis and Lemma \ref{lemma:boundend_deterioration_complete}. \fin

It is interesting to note that, in the particular case that $G$ is contractive, Theorem~\ref{th:AAconvergence_theorem} proves that, at least locally,
the stabilized version of AA (Algorithm~\ref{alg:ATM_2}) could improve the rate of convergence of the fixed point map
$\mathbf{s}_{j+1}=G(\mathbf{s}_j)$ since the linear convergence parameter $q$ in Theorem~\ref{th:AAconvergence_theorem}
can be chosen smaller than the contraction factor of $G$ (see also \cite{Evans,C}). Observe, moreover, that
 if the inequality $\|-\beta I-J^{-1}\|\leq \delta$ could not be fulfilled, we can
consider the \textit{preconditioned} non linear function $\tilde{F}= P^{-1}F$
where $P$ is some approximation of $J=F(\mathbf{s}^*)$, and we obtain in this way
$\|-\beta I-J^{-1}P\|\leq \delta$.

Finally, let us observe that, as customary in the quasi-Newton literature, we can
improve  the global convergence properties of the {AA} procedure by
 introducing  a step-length parameter $\alpha_j$
 and {transforming} the sequence generated by  Algorithm~\ref{alg:ATM_2} into the sequence
\begin{equation*}
	\mathbf{s}_{j+1}=\mathbf{s}_j-\alpha_j H_j^{(\beta)} \mathbf{f}_j .
\end{equation*}

\subsection{Connections between stabilized AA and regularized ATM} \label{sec:stabilization}
As already pointed out in the previous section, from a theoretical point of view, the stabilization procedure introduced in
Algorithm~\ref{alg:ATM_2}, in order to ensure the convergence, aims to detect a subset of the vectors in $\Delta
F_{j-m_j}^{(m_j)}$ that are \textit{sufficiently linearly  independent}: the proposed stabilization procedure in
Algorithm~\ref{alg:ATM_2} (Lines \ref{alg:lineli}
  -\ref{alg:linelif}) can be interpreted simply as a Gram-Schmidt
procedure with threshold, i.e., the residual difference $\Delta \mathbf{f}_{d}$ is discarded if it is \textit{close} to a vector
linearly
  dependent from the previously computed residual differences. The
above observation naturally links the stabilization procedure with Rank-Revealing QR factorizations \cite{chanRRQR,guSRRQR}.  We find
this issue particularly interesting and deserving further investigation. Here, we prefer
to adopt a regularization point of view, as in
\cite{aspre,ander1,B,C}, to motivate the introduction of the regularization parameter $\lambda$ in the Anderson-Type
Mixing methods as we did at the beginning of Section \ref{sec:anderson_type}. To this end,
let us consider the ATM obtained by
\eqref{eq:dopo34} with ${\mathbf c}_{j} = {\mathbf f}_{j}$ and $M=I$. As already pointed out, when $\lambda=0$ it
coincides with the classical AA but, when $\lambda \neq 0$, the method obtained can be viewed as a Regularized Anderson Acceleration (RAA).

In this setting, we interpret the magnitude of the singular values of the matrix $\Delta
	F_{j-m_j}^{(m_j)}$  as a measure of the linear independence of its columns:
the presence of linearly dependent vectors in $\Delta F_{j-m_j}^{(m_j)}$ is highlighted by the presence of very small singular values.
Let us consider now the SVD decomposition  $\Delta F_{j-m_j}^{(m_j)}=U\Sigma V^T$.
We add a regularization parameter $\lambda$ to the matrix $\Sigma$ and, we set
$\Delta \widetilde{F}_{j-m_j}^{(m_j)} = U\sqrt{\Sigma^2+\lambda I }V^T $.
By direct computation, it is possible to show that \eqref{eq:dopo34} can be written as
\begin{equation*}
\boldsymbol{\theta}^{(j)}_{I,\lambda}=((\Delta \widetilde{F}_{j-m_j}^{(m_j)} )^T \Delta \widetilde{F}_{j-m_j}^{(m_j)} )^{-1}(\Delta
F_{j-m_j}^{(m_j)})^T{\mathbf f}_{j}.
\end{equation*}
The statement regarding the linear independence of the columns of the matrix $\Delta \widetilde{F}_{j-m_j}^{(m_j)}$
can be obtained by observing that all its singular values are bounded from below by $\sqrt{\lambda}$.
We consider the above argument as an explanation of the fact that the introduction of a
regularization parameter in  the AA method (and, in general, in all the ATMs)  could achieve numerically the same task of the
stabilization procedure of Algorithm~\ref{alg:ATM_2}.
Adopting a quasi-Newton point of view, it is important to observe that using Formula \eqref{eq:general_theta_AA_UPDATE}
with $\boldsymbol{\theta}^{(j)}=\boldsymbol{\theta}^{(j)}_{I,\lambda}$, the ATM update (see Line~\ref{algline:AAupdate} in
Algorithm~\ref{ATMtable_1}) can be written as
\begin{equation*}
\mathbf{s}_{j+1}=\mathbf{s}_{j}-\widetilde{H}_j^{(\beta)}\mathbf{f}_j,
\end{equation*}
with
\begin{equation}\label{eq:anderson_type2}
	\widetilde{H}_j^{(\beta)}=  - \beta I +(\Delta S_{j-m_j}^{(m_j)}+ \beta \Delta F_{j-m_j}^{(m_j)})(( \Delta \widetilde{F}_{j-m_j}^{(m_j)})^T \Delta \widetilde{F}_{j-m_j}^{(m_j)})^{-1}(\Delta F_{j-m_j}^{(m_j)})^T.
	\end{equation}
The quasi-Newton matrices defined in \eqref{eq:anderson_type2} satisfy {only} an \textit{approximated multisecant condition},
	namely \[
	\widetilde{H}_j^{(\beta)}
	\Delta F_{j-m_j}^{(m_j)} =\Delta S_{j-m_j}^{(m_j)}+ \beta (\Delta F_{j-m_j}^{(m_j)} (( \Delta \widetilde{F}_{j-m_j}^{(m_j)})^T \Delta \widetilde{F}_{j-m_j}^{(m_j)})^{-1} \Delta {F_{j-m_j}^{(m_j)}}^T\Delta F_{j-m_j}^{(m_j)}- \Delta F_{j-m_j}^{(m_j)} ),
	\]
which represents a noteworthy difference with the multisecant conditions satisfied by the quasi-Newton matrices used in the classical AA
and in its stabilized version (see Lemma~\ref{lem:Multisecant_stabilized}).

\section{Numerical results} \label{sec:numerical_res}
In this Section, we investigate the numerical behavior of some of the methods studied in the previous sections for different test problems.

\subsection{Details on the methods and their implementations}
We select a subset of the methods presented in the previous sections with the main aim to compare their numerical
  performance (with a focus on the rate of convergence), and to prove that the acceleration performance they deliver
  behave consistently. Our choices are, among other things, driven by the fact that all the acceleration methods considered share the
  same order of complexity (linear in the dimension
  of the problem) per acceleration step.  A comprehensive detailed
  numerical study and the relative implementations of all the methods
  described  in the previous sections is out of the scope of this
  work and is postponed to future works.   Table
  \ref{tab:methods} summarizes the methods we consider in our
  numerical experiments. In the first column we report the name and
  the relative abbreviation for the particular acceleration scheme we
  consider. In the second column we report the reference equations of
  the acceleration scheme and, for the sake of completeness, in the
  third column we report the strategy type of the considered
  acceleration: \textit{Restarted Method} (RM) or \textit{Continuous-Updating}
  (CU).  Finally, in the last column, we report the details
  concerning the choice of the regularization parameter: in the \textit{Grid
  Search} (GS) approach the regularization parameter $\overline{\lambda}$ is chosen, as proposed in
  \cite{aspre}, as the parameter which achieves the smallest
  fixed point residual;  {the interval $[10^{-12},1]$ is discretized logaritmically into $7$ values of $\lambda$,  (for more details see Algorithm \ref{RMtable_reg_selection} which is a modification of Algorithm \ref{RMtable}) among which, one of them, $\overline{\lambda}$, is selected}. For the sake of completeness, let us recall that,
  	also in this new algorithm,  $\ell_k=k+2$ if we use
  	\eqref{st} or \eqref{eq:beta_expression}, and $\ell_k=2k+1$ if we use
  	\eqref{eq:topological_extrapolated}
  	or \eqref{eq:topological_extrapolated_beta_new}
  For the \textit{Generalized Cross Validation} (GCV) approach,
  which is a natural approach for regularizing ill-posed
  regression-like problems, we refer the interested reader to
  \cite{GCVgolub}.

\begin{table}[ht!]
  \footnotesize
  \begin{center}
    \begin{tabular}{|l|l|l|l|}
      \hline
      \textbf{Name}  &  \textbf{Ref. Eq.}  &  \textbf{Type}               & \textbf{Choice of $\lambda$}        \\ \hline	
      Singular Value Decomposition Acceleration ({SVDA})	        &\eqref{eq:svd} & RM & $\lambda=0$  \\ \hline
      Regularized Nonlinear Acceleration ({RNA})	                &\eqref{st0} & RM &GS (Alg. \ref{RMtable_reg_selection})  \\ \hline
      Regularized Reduced Rank Extrapolation ({RRRE})	            &\eqref{eq:beta_expression0} & RM &GCV \cite{GCVgolub}  \\ \hline
      Regularized Topological Shanks Acceleration ({RTSA})	     &  \eqref{eq:solstar}   & RM &GS (Alg. \ref{RMtable_reg_selection}) \\ \hline

      Anderson Acceleration ({AA}) with $0<\beta \leq 1$	                                &\eqref{eq:quasi_newton_variable_beta} & CU & $\lambda=0$  \\ \hline
      Regularized Anderson Acceleration ({RAA})                   & \eqref{eq:anderson_type2}& CU &GCV \cite{GCVgolub}\\ \hline
    \end{tabular}
  \end{center}	
\caption{Methods tested. \label{tab:methods}}
\end{table}

Let us point out that, to the best of our knowledge, among the methods presented in Table~\ref{tab:methods}, RTSA and RRRE/RAA with the regularization parameter chosen using the GCV are new approaches
introduced in this work. Instead, for the other methods, we refer in particular to \cite{sidisvd} for the SVDA
(which is called SVD-MPE in the original paper) and to \cite{aspre} for the RNA.

Finally, we mention  that in all the numerical experiments we used $M=I$ and that,
in the SVDA approach, we use as extrapolated term
	${\mathbf t}_n^{(k+1)}=S_{n+1}^{(k+1)}\boldsymbol{\alpha}$ where $\boldsymbol{\alpha}$ is the normalized singular vector corresponding to the smallest singular value of $\Delta S_n^{(k+1)}$ (see equation \eqref{eq:svd}).

\begin{algorithm}[ht]
	\caption{The Restarted Method (RM) with grid-search (GS).}\label{RMtable_reg_selection}
	\DontPrintSemicolon
	\KwIn{Choose  $M$, $k$, $\lambda_{\min}$, $\lambda_{\max}$, $n$, and ${\mathbf x}_0 \in \mathbb{R}^p$.}
	\For{$j=0,1,\ldots$}{Set ${\mathbf s}_0= {\mathbf x}_j$\;
		\For{$i=1,\ldots,\ell_k-1$ (\it basic or inner iterations)}{
			Compute ${\mathbf s}_i=G({\mathbf s}_{i-1})$ \;}
		Choose $\lambda_0, \ldots, \lambda_{n-1} \in [\lambda_{\min},\lambda_{\max}]$\;
		\For{$i = 0, \ldots, n-1$}{
			Compute $\mathbf{t}_{0,\lambda_i}^{(\ell_k-1)}$ using \eqref{st}~or~\eqref{eq:beta_expression}~or~\eqref{eq:topological_extrapolated}
			~or~\eqref{eq:topological_extrapolated_beta_new}  \;}
			 {${\overline{\lambda}}= \underset{\lambda_i \in \{\lambda_0,\dots,\lambda_{n-1}\}}{\arg \min} \|G(\mathbf{t}_{0,\lambda_i}^{(\ell_k-1)})-\mathbf{t}_{0,\lambda_i}^{(\ell_k-1)}\| $}\;
			Set $\mathbf{x}_{j+1}=\mathbf {t}_{0,\overline{\lambda}}^{(\ell_k-1)}$\;
}
\end{algorithm}

All the numerical experiments are performed on a laptop running Linux
with 16Gb memory and CPU Intel\textsuperscript{\textregistered}
Core\texttrademark\ i7-4510U with a clock speed of 2.00GHz. The code is written
and executed in Python. For the discretization of the PDE's we used
Fenics \cite{fenics} and, for the GCV choice of the regularization
parameter, we used the Scikit-learn package \cite{scikit}.
 {Throughout the experiments, to show and test the robustness of the different extrapolation approaches, we base all our extrapolation schemes on {$7$} previous iterates, i.e., $\ell_k=7$ in Algorithm \ref{RMtable_reg_selection} or $m=7$ in Algorithm \ref{ATMtable_1} }

\subsection{PageRank}
The aim of this first numerical example is to highlight the benefits
of introducing regularization strategies in Shanks-based
extrapolation methods. In particular, in this section, we consider
the PageRank problem (see \cite{eldendata}), i.e., the problem of
computing the Perron eigenvector of the matrix
\begin{equation*}
	G=\alpha S+\frac{(1-\alpha)}{n}\mathbf{e}\mathbf{e}^T, \;\;\alpha\in (0,1),
\end{equation*}
where $S$ is a nonnegative column stochastic matrix.  For the solution
of this problem, we consider the Power Method, i.e.,
$\mathbf{u}_{k+1}=G(\mathbf{u}_k)$ where $\mathbf{u}_0$ is a
nonnegative stochastic vector, which is known to be a linear fixed
point iteration globally convergent with a rate of convergence of
$O(\alpha^{k})$ \cite{eldendata}. As the previous convergence bound
confirms, the rate of convergence of the Power Method for the PageRank
computation becomes slower as $\alpha$ approaches $1$, but this is
usually the case of interest in applications \cite{eldendata}. In this experiment we use as stopping criterion $\|G(\mathbf{u}_k)-\mathbf{u_k}\|< 10^{-7}$.

\begin{figure}[th!]
	\centering
	\includegraphics[width=.45\textwidth]{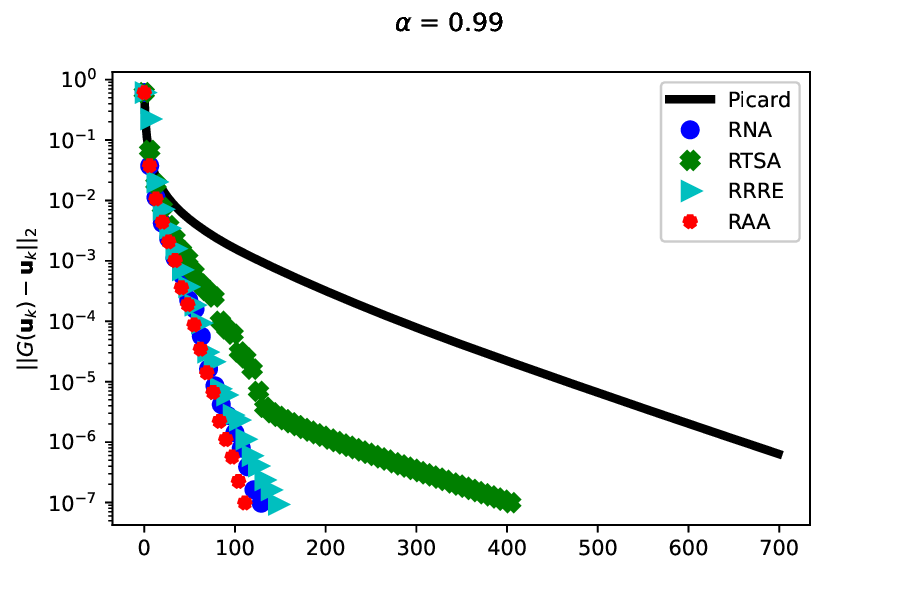}
	\includegraphics[width=.45\textwidth]{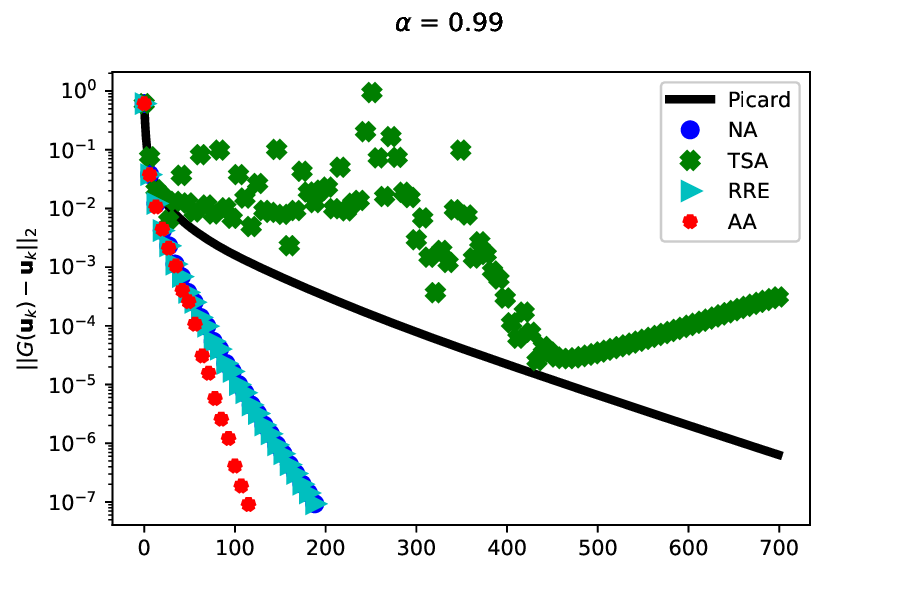}\\
	\includegraphics[width=.45\textwidth]{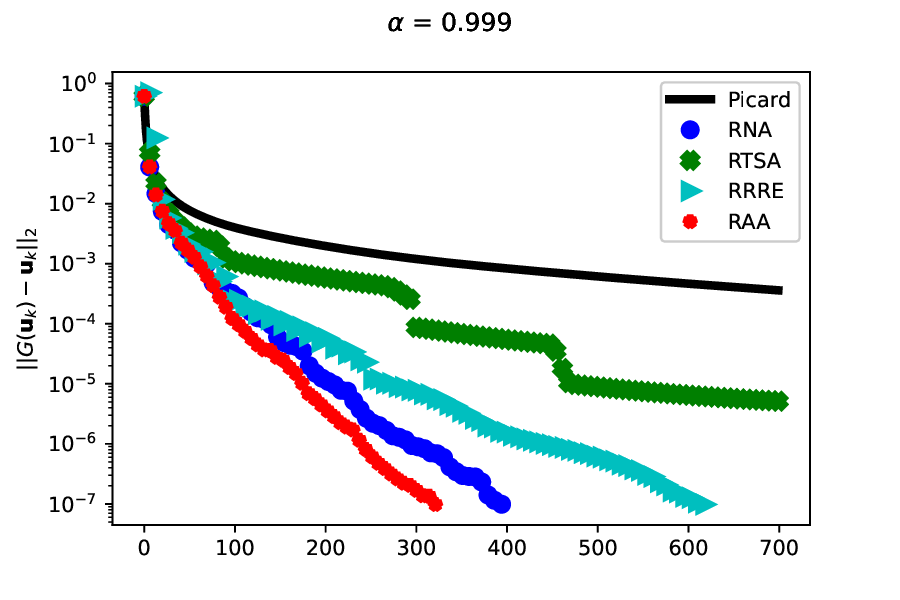}
	\includegraphics[width=.45\textwidth]{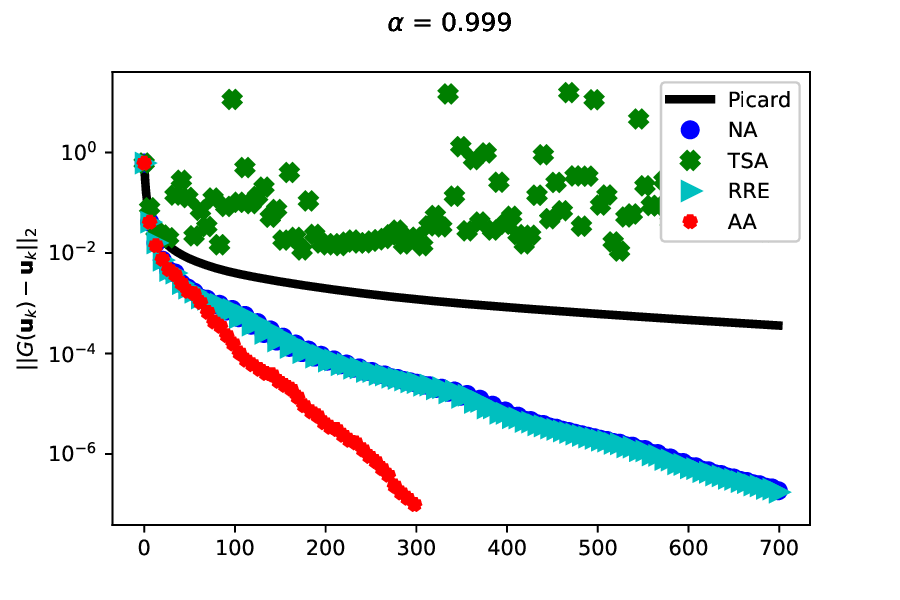}\\
	\caption{PageRank Problem} \label{fig:PRproblem}
\end{figure}

In the left panel of Figure \ref{fig:PRproblem}, we report the acceleration performance of the regularized versions
of the methods considered when compared to the non regularized ones (in the right panel), for the computation of the
PageRank vector of the matrix \texttt{amazon-0202} from \cite{davismatrix} (which has been suitably modified in
order to be stochastic and Dangling-Nodes free \cite{eldendata}). Recall that the sequence generated by
the Power Method belongs to the Shanks kernel and hence,
at least theoretically,  all the extrapolation strategies should be equivalent and should work consistently without
any requirement of regularization. Nevertheless, as Figure \ref{fig:PRproblem} clearly shows, the introduction
of a regularization strategy improves the robustness of the extrapolation procedures permitting, for the restarted
extrapolation methods (namely RNA, RTSA, RRRE), to obtain a more  effective acceleration performance
across different choices of the parameter $\alpha$. Observe also that, in this case, the introduction
of a regularization procedure in the {AA} scheme ({RAA}) does not sensibly improves the acceleration performance.

\subsection{Nonlinear Poisson problems}
In this Section, we consider the solution of the nonlinear PDE (see equation \eqref{eq:NLPDE_problem})

\begin{equation} \label{eq:NLPDE_problem}
\begin{split}
- \nabla (q(u) \nabla u) + g(u)  +u_x & = f \hbox{ in } \mathcal{D}=[0,1] \times [0,1], \\
u & = v \hbox{ on } \partial \mathcal{D}.
\end{split}
\end{equation}
We use a $1/64$ uniform triangular mesh of $\Omega=[0,1]^2$ with a $(P_2)$ discretization \cite{fenics} that provides a total of $16,641$ degrees of freedom.
In particular, we consider the following choices of the functions

\begin{itemize}
	\item $q(u)=1+u^2$ or $q(u)=1+u^4$, $g(u) =  0$ and $f$ such that the exact solution of \eqref{eq:NLPDE_problem}
is given by $\overline{u}=   \exp(-2x)\sin(3\pi y)$ and $v=\overline{u} \hbox{ on } \partial \Omega$.
We refer to these choices as the \textit{Nonlinear Poisson Problem};
	\item $q(u)=   1$, $g(u)= \lambda e^{u}$ with $\lambda=1$ or $\lambda =-1$, $f =   0$
and $v=   0 \hbox{ on } \partial \Omega$. We refer to these choices as the \textit{Bratu Problem} \cite{hjbbratu}.
\end{itemize}

\begin{figure}[th!]
	\centering
	\includegraphics[width=.45\textwidth]{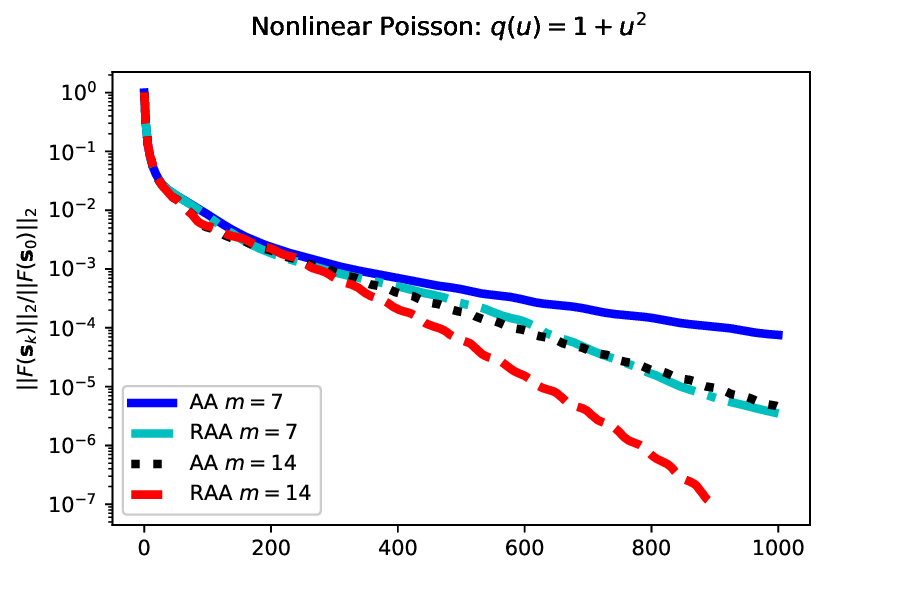}
	\includegraphics[width=.45\textwidth]{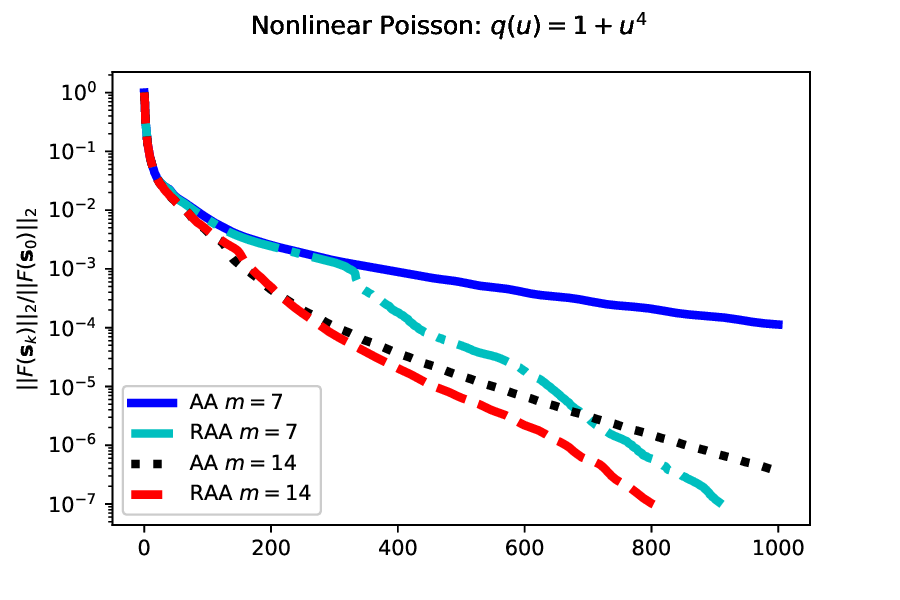}\\
	\caption{Nonlinear Poisson Problem.} \label{fig:NLP}
\end{figure}

\begin{figure}[th!]
	\centering
	\includegraphics[width=.45\textwidth]{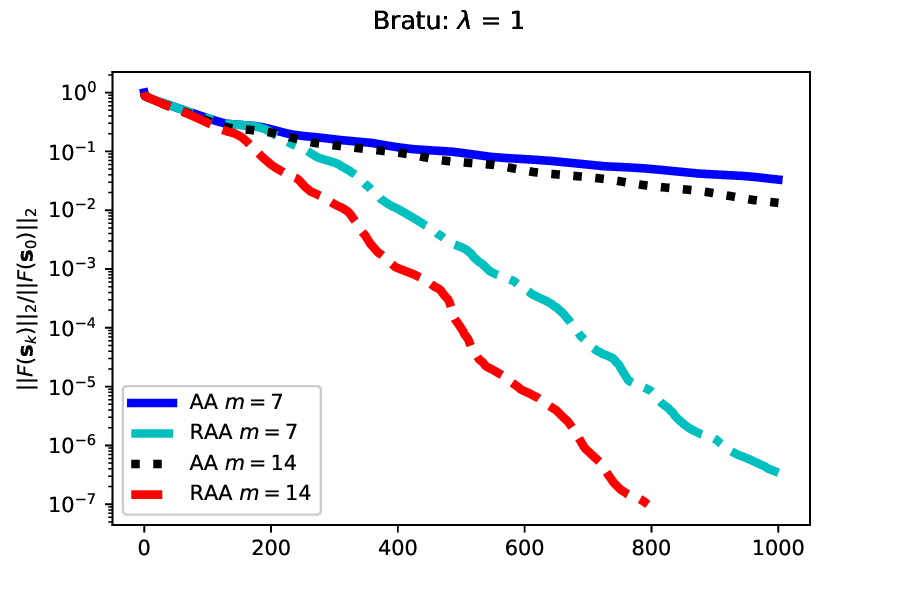}
	\includegraphics[width=.45\textwidth]{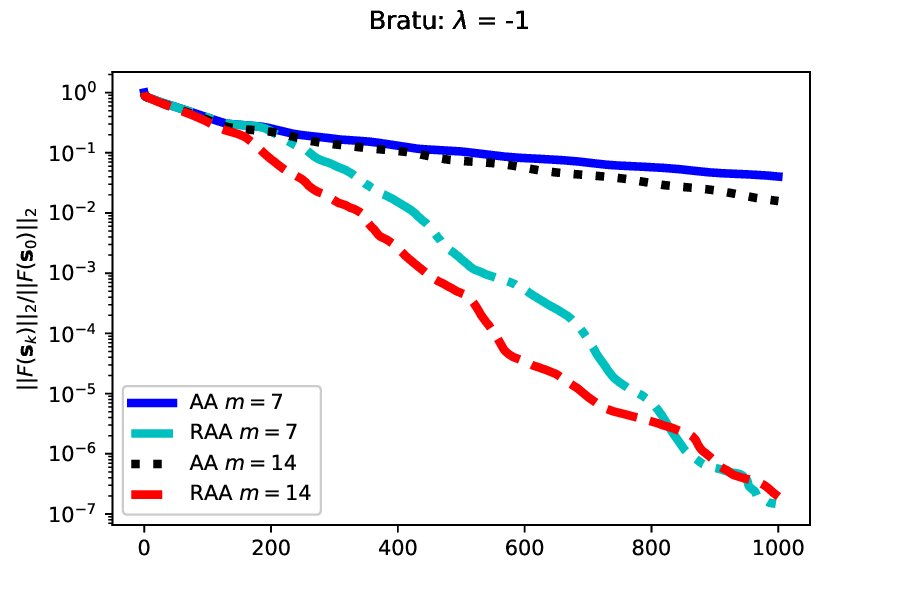}\\
	\caption{ Bratu Problem.} \label{fig:Bratu}
\end{figure}

After the discretization of \eqref{eq:NLPDE_problem}, the
corresponding problems can be written as the solution of
$F(\mathbf{s})=0$, i.e., as the solution of a linear system of
equations. We assume that the derivative of $F$ are not readily
available or that a sufficiently accurate initial guess is not at our
disposal in order to apply Newton's method. In this experiment we use as stopping criterion $\|F(\mathbf{u}_k)\|< 10^{-7}$. Figures~\ref{fig:NLP} and \ref{fig:Bratu}, show  the acceleration
performance of {AA} when compared to its regularized version
{RAA} (these problems are not well scaled and a good choice for
the mixing parameter was $\beta=0.1$) for the problems previously
discussed. The Figures clearly show that the introduction of the
regularization strategy, in these cases characterized by a
\textit{higher nonlinearity} than for the  PageRank example, leads to a better
robustness of the schemes with respect to the choice of the memory
parameter $m$. In particular, the introduction of the
regularization procedure permits to have a satisfactory
rate of convergence independently from the value $m$.
We point out that, interestingly enough, the need for a stabilization
procedure needed from the theoretical point of view to prove the
convergence of the {AA} scheme (see Algorithm \ref{alg:ATM_2}),
is echoed by the experimental observation that increasing
$m$ could result in a lost of efficiency for the AA
scheme (see Figure \ref{fig:Bratu}). The introduction of a
regularization procedure mitigates such a drawback.

\begin{figure}[th!]
	\centering
	\includegraphics[width=.45\textwidth]{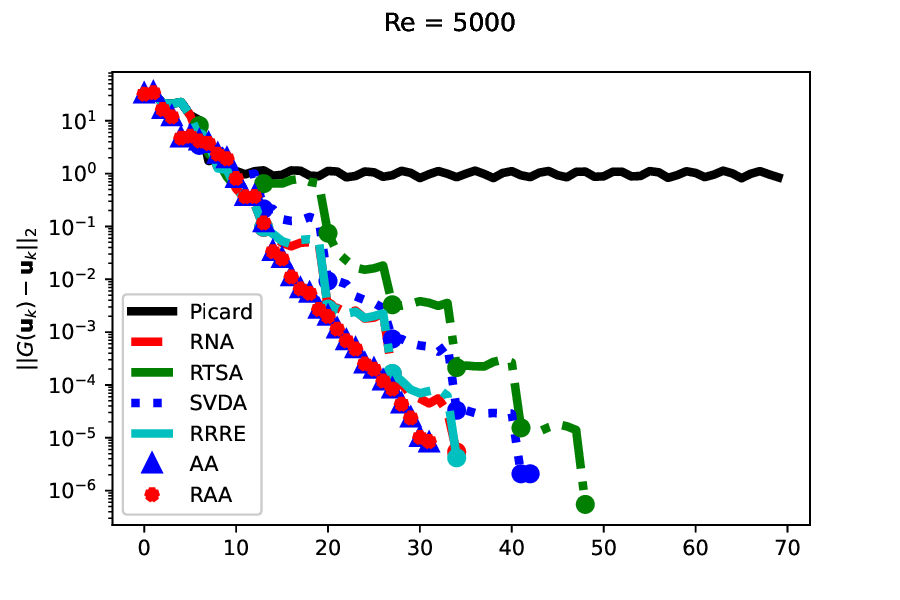}
	\includegraphics[width=.45\textwidth]{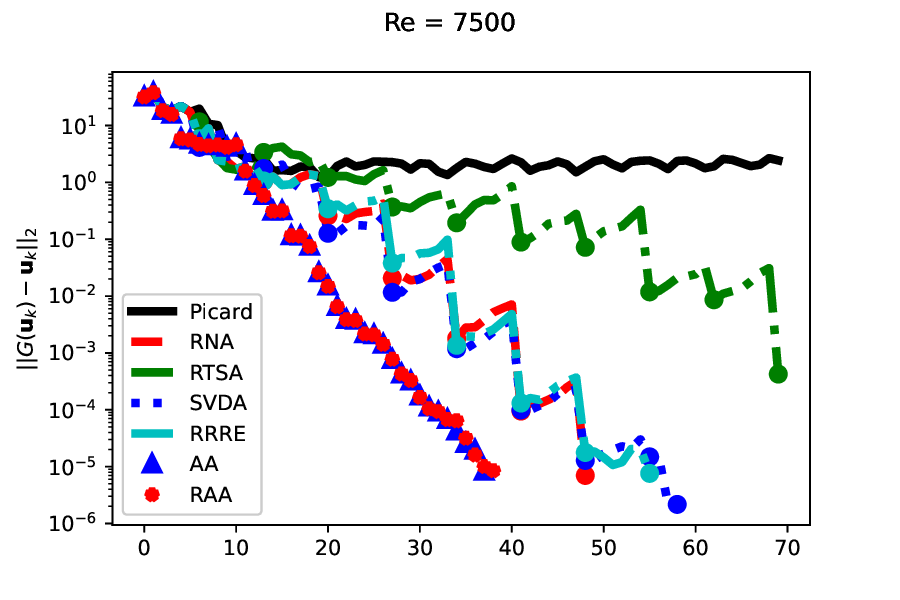}\\
	\includegraphics[width=.45\textwidth]{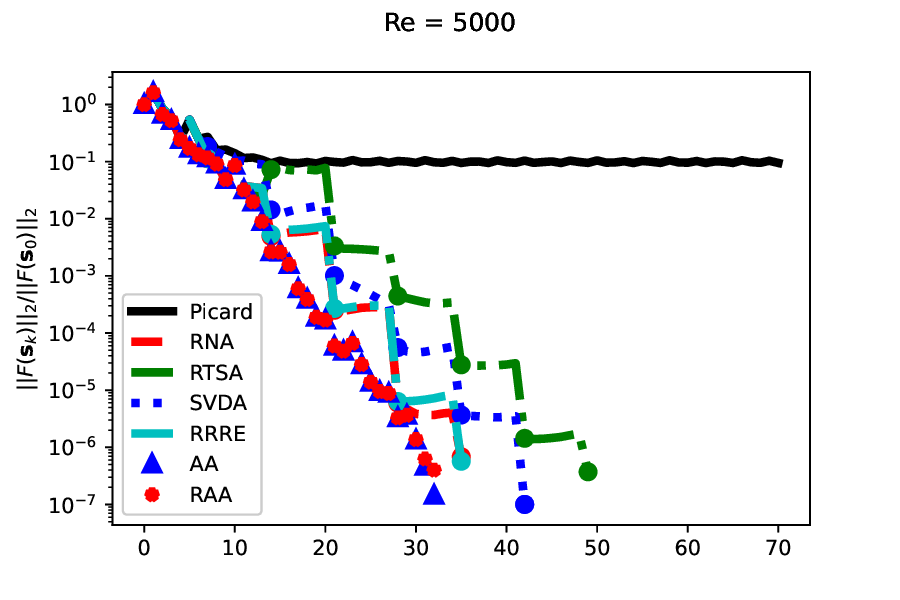}
	\includegraphics[width=.45\textwidth]{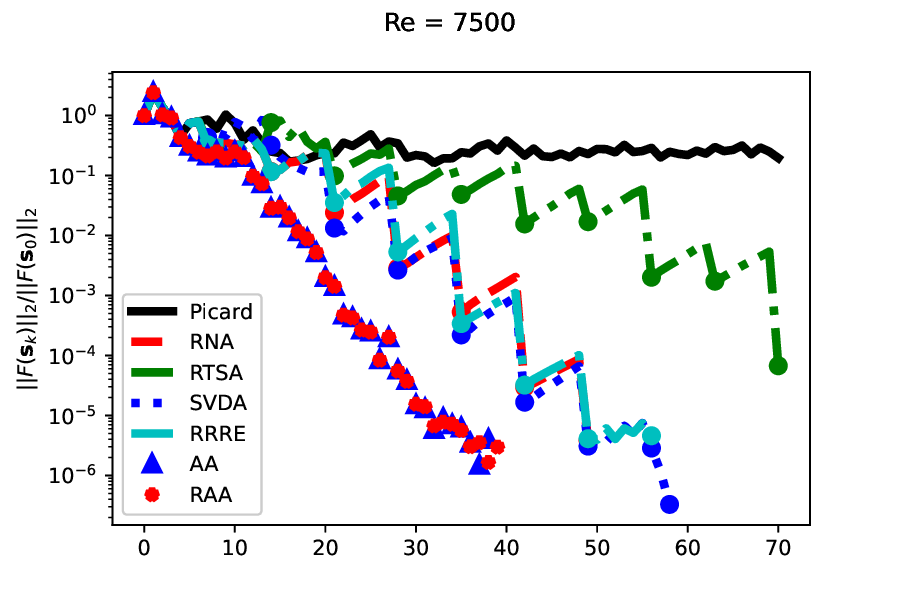}\\
	\caption{ Lid driven Problem: Acceleration Performance.} \label{fig:RMNSE}
\end{figure}

\begin{figure}[th!]
	\centering
	\includegraphics[width=.45\textwidth]{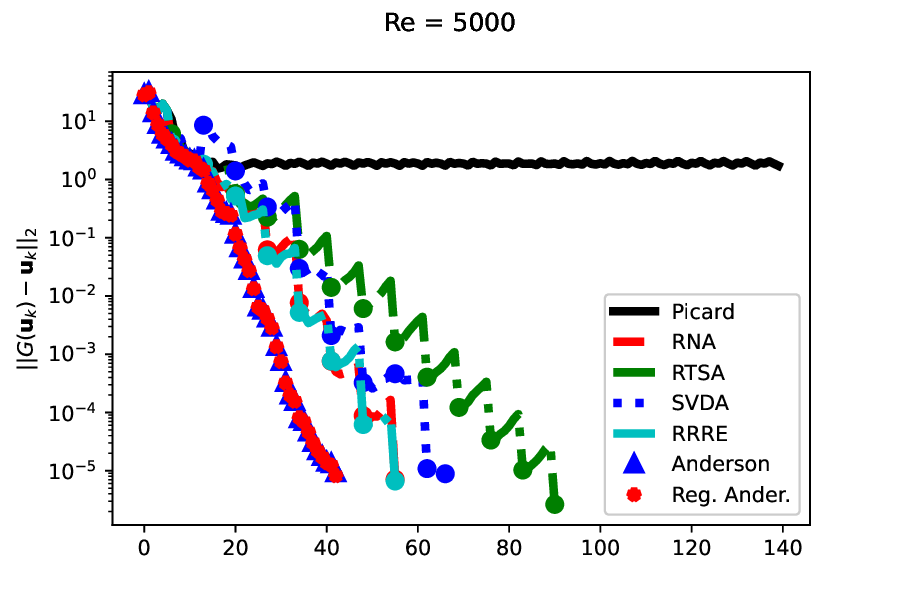}
	\includegraphics[width=.45\textwidth]{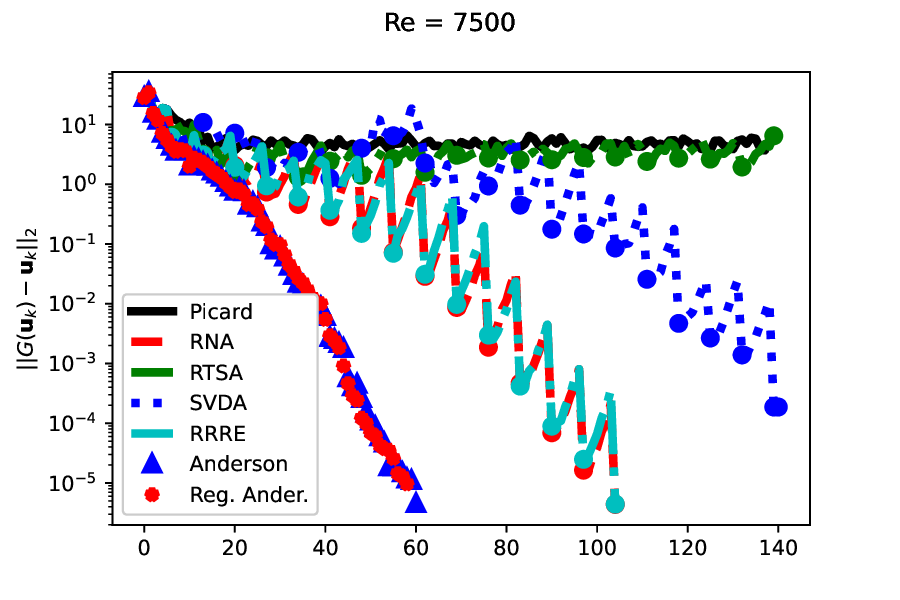}\\
	\includegraphics[width=.45\textwidth]{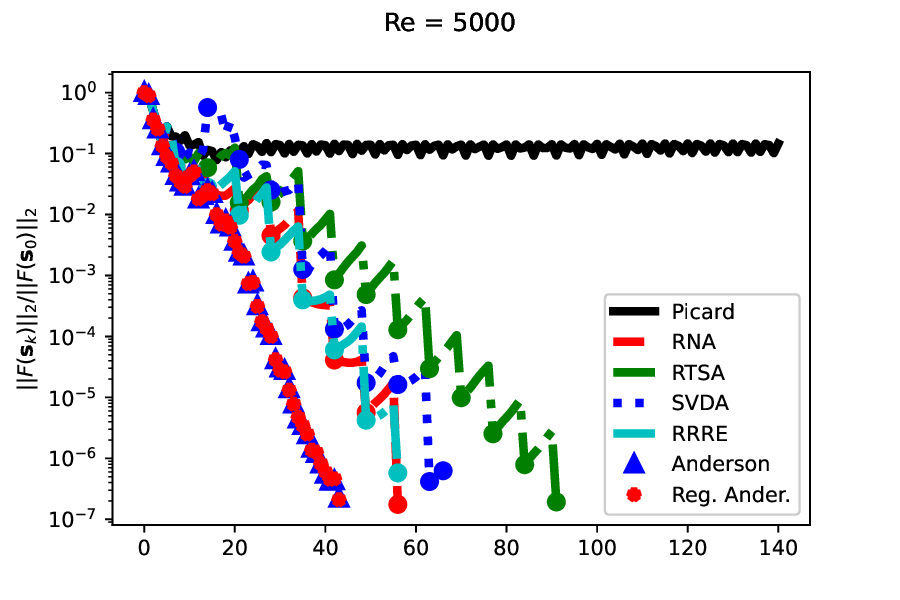}
	\includegraphics[width=.45\textwidth]{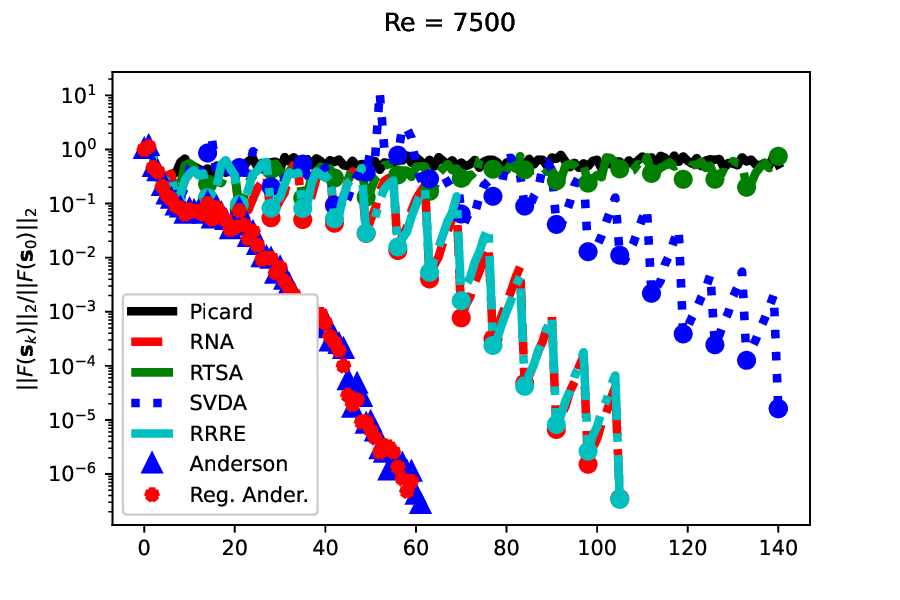}\\
	\caption{ Lid driven Problem: Acceleration Performance.} \label{fig:RMNSEdeep}
\end{figure}
\subsubsection{Navier-Stokes equation}
In this Section, we compare the numerical performance of the different restarted
extrapolation approaches on the incompressible Navier-Stokes
Equation (NSE)
\begin{align*}
u \cdot \nabla u+ \nabla p- \nu \Delta u= f,\\
\nabla \cdot u =0, \\
 u|_{ \partial \Omega}=g,
\end{align*}
where $\nu$ is the kinematic viscosity, $f$ is the forcing, $u$ and
$p$ represent velocity and pressure and $\Omega$ is a given domain in
$\mathbb{R}^2$.
Following \cite{pollock1},
we consider a Picard iteration {(equations \eqref{eq:picardNSE2})} to solve the problem.
The iteration, which is commonly used for its stability and global convergence properties, takes  the form

\begin{eqnarray}
u_k \cdot \nabla u_{k+1}+ \nabla p_{k+1}- \nu \Delta u_{k+1}= f,&&
\nonumber \\
\nabla \cdot u_{k+1} =0, &&\label{eq:picardNSE2}\\
u_{k+1}|_{\partial \Omega}=g . &&
\nonumber
\end{eqnarray}

The above scheme is written in the fixed point form $u_{k+1}=G(u_k)$,
where $G$ denotes the solution operator for the linearization
\eqref{eq:picardNSE2}. To be specific, we
consider the 2D lid driven cavity ($\Omega=(0,1)^2$) and a ``deep'' lid
driven cavity with ($\Omega=(0,1)\times (0,3)$). No slip ($u=0$)
boundary conditions are imposed on the sides and the bottom, and the
Dirichlet boundary condition $u(x,1)=(1,0)^T$ is imposed on the top
to enforce the ``moving lid'' condition. There is no forcing $(f=0)$
and the kinematic viscosity ($\nu =  Re^{-1}$) is considered at
benchmark values $Re=5000,\;7500$. We discretize with $(P_2,P_1)$
Taylor Hood elements. In the case $\Omega=(0,1)^2$ we use a
$\frac{1}{64}$ uniform triangular mesh that provides a $37,507$ total
degrees of freedom and in the case $\Omega=(0,1)\times (0,3)$ we use a
$\frac{1}{40} \times \frac{1}{120} $ mesh that provides $87,203$ total
degrees of freedom.  Similarly to the results presented in
  \cite{pollock1}, our experiments confirm that  Newton's method
  starting with a zero initial guess, never converges.  In this experiment we use as stopping criterion $\|G(\mathbf{u}_k)-\mathbf{u_k}\|< 10^{-5}$. Figures
  \ref{fig:RMNSE} and \ref{fig:RMNSEdeep} show the acceleration
  performance of the methods described in Table \ref{tab:methods} for
  the solution of the steady NSE.  The best performer in terms of
  achieved acceleration is {AA} and the introduction of a
  regularization procedure in this scheme ({RAA}) seems not to
  have a relevant impact on the rate of convergence. This is
  probably due to the fact that the fixed point iteration we are
  considering generates a sequence that is \textit{close} to being a linear
  sequence and, as in the PageRank case, regularization of
  the {AA} scheme does not seem to have a great influence. Concerning
  the restarted regularized methods, we should notice that the
  {RTSA} is not able to achieve an acceleration performance in the
  \textit{Deep} case for $Re=7500$. Finally, let us highlight the
  particularly interesting performance of the {SVDA} approach:
  this approach does not require the computation of any regularization
  parameter and only one SVD decomposition every  {$\ell_{k}-1$} fixed point
  iterations is needed, whereas {AA} requires the solution of a
  least square solution per step, and all the regularized methods which
  use the RM approach require the selection of a regularization
  parameter. The non-regularized versions of the methods using the RM
  strategy, as in the PageRank case, exhibited a worst performance and
  are not reported for this reason.

\section{Conclusions}

In this work, we presented a unified framework for Shanks-based
transformations.  If, on one hand, the introduction of this
framework allowed us to link apparently different
extrapolation/acceleration techniques with Shanks-based
transformations, on the other hand, it allowed us to introduce suitable
generalizations able to numerically outperform the existing ones, as
highlighted in the preliminary numerical results presented,
  especially on problems characterized by a high degree of
  nonlinearity. To conclude, we note  that the
  highlighted connection between the Shanks-based transformations and the
  quasi-Newton methods and Anderson Acceleration,   shed
   light into some of its theoretical and numerical behaviors,
   furthering our knowledge of  the powerful, but poorly understood, Anderson acceleration
   \cite{kelley2018numerical}.

\section*{Acknowledgments:}
We would like to thank the reviewers for their very careful reading of our paper, and for their constructive comments.

C.B. acknowledges support from the Labex CEMPI (ANR-11-LABX-0007-01).
S.C. and M.R.-Z. are a members of the INdAM Research group GNCS. The work of S.C. was partially supported by the GNCS –
INdAM project \textit{``Efficient Methods for large scale problems with applications to data analysis
  and preconditioning''} and from Dept. of Computer Science \& Engineering, University of
  Minnesota, \textit{Project No. UMF0002384}. The work of M.R.-Z. was partially supported by the University of Padua, Project No.
DOR 1903575/19 \textit{Numerical Linear Algebra and Extrapolation methods with applications}. The work of Y.S. was supported in part by NSF grant DMS-1912048.


\begin{thebibliography}{99}
\bibitem{fenics}
M. Aln\ae s, J. Blechta, J. Hake, A. Johansson, B. Kehlet, A. Logg, C. Richardson, J Ring, M.E. Rognes, G.N. Wells,
The FEniCS project version 1.5, Archive of Numerical Software, 3(100) (2015) 9-23.

\bibitem{ander}
D.G. Anderson, Iterative procedures for nonlinear integral equations, J. Assoc. Comput.
Mach., 12 (1965) 547-560.
\bibitem{ander1}
D.G. Anderson, Comments on ``Anderson acceleration, mixing and
extrapolation'', Numer.  Algorithms, 80 (2019) 135–234.

\bibitem{A}
A.S. Banerjee, P. Suryanarayana, J.E. Pask,   Periodic Pulay method for robust and efficient convergence acceleration
of self-consistent field iterations. Chemical Physics Letters, 647 (2016) 31-35.

\bibitem{cb1}
C. Brezinski, Application de l'$\varepsilon$-algorithme \`a la r\'esolution des
syst\`emes non lin\'eaires, C. R. Acad. Sci. Paris, 271A (1970) 1174-1177.

\bibitem{cbth}
C. Brezinski, {\it M\'{e}thodes d'Acc\'{e}l\'{e}ration de la Convergence en Analyse Num\'{e}rique},
Th\`{e}se de Doctorat d'\'{E}tat, Universit\'{e} Scientifique et M\'{e}dicale de Grenoble, 1971.
https://tel.archives-ouvertes.fr/tel-00282774
\bibitem{vec}
C. Brezinski, Some results in the theory of the vector $\varepsilon$-algorithm,
Linear Algebra Appl., 8 (1974) 77-86.
\bibitem{etopo}
C. Brezinski, G\'en\'eralisation de la transformation de Shanks, de la table
de Pad\'e et de l'$\varepsilon$-algorithme, Calcolo, 12 (1975) 317-360.
\bibitem{birk}
C. Brezinski, {\it Pad\'e-Type Approximation and General Orthogonal Polynomials},
ISNM, vol. 50, Birkh\"auser-Verlag, Basel, 1980.

\bibitem{sch}
C. Brezinski, Other manifestations of the Schur complement, Linear Algebra
Appl., 111 (1988) 231-247.

\bibitem{proj}
C. Brezinski, {\it Projection Methods for Systems of Equations}, Elsevier, Amsterdam, 1997.

\bibitem{cbmrz}
C. Brezinski, M. Redivo-Zaglia, {\it Extrapolation Methods. Theory and
Practice}, North-Holland, Amsterdam, 1991.


\bibitem{simply}
C. Brezinski, M. Redivo-Zaglia, The simplified topological
$\varepsilon$-algorithms for accelerating sequences in a vector space, SIAM J. Sci.
Comput., 36 (2014) A2227-A2247.

\bibitem{soft}
C. Brezinski, M. Redivo-Zaglia, The simplified topological $\varepsilon$-algorithms: software and applications,
Numer. Algorithms, 74 (2017) 1237-1260.

\bibitem{genesis}
C. Brezinski, M. Redivo-Zaglia, The genesis and early developments of Aitken’s process, Shanks’ transformation,
 the $\varepsilon$-algorithm, and related fixed point methods, Numer. Algorithms, 80(1) (2019) 11-33.

\bibitem{era}
C. Brezinski,  M. Redivo-Zaglia, {\it Extrapolation and Rational Approximation. The Works of the Main Contributors},
Springer Nature, Cham, Switzerland, 2020.

\bibitem{ext}
C. Brezinski,  M. Redivo-Zaglia, Extrapolation and prediction of sequences in a vector space, submitted.


\bibitem{brzs}
C. Brezinski, M. Redivo-Zaglia, Y. Saad, Shanks sequence transformations and Anderson acceleration,
SIAM Rev., 60(3) (2018) 646-669.

\bibitem{badBroyden}
C.G. Broyden, A class of methods for solving nonlinear simultaneous
equations, Math. Comp., 19 (1965) 577-593.
\bibitem{mpe}
S. Cabay, L.W. Jackson, A polynomial extrapolation method for finding limits and antilimits of vector sequences,
SIAM J. Numer. Anal., 13 (1976) 734-752.
\bibitem{chanRRQR}
T.F. Chan, Rank revealing QR factorizations, Linear Algebra Appl., 88 (1987) 67-82.

\bibitem{cipollamulti}
S. Cipolla, M. Redivo-Zaglia, F. Tudisco, Extrapolation methods for fixed-point multilinear PageRank computations, Numer. Linear Algebra. Appl.,  27 (2020)  e2280.
\bibitem{cipollashifted}
S. Cipolla, M. Redivo-Zaglia, F. Tudisco, Shifted and extrapolated power methods for tensor $\ell^p$-eigenpairs, Electron. Trans. Numer. Anal., 53 (2020) 1-27.


\bibitem{davismatrix}
T.A. Davis, Y. Hu, The University of Florida sparse matrix collection, ACM Trans. Math. Software, 38(1) (2011) 1-25.


\bibitem{dela}
J.P. Delahaye, {\it Sequence Transformations}, Springer-Verlag, Berlin, 1988.
\bibitem{bgbj}
J.P. Delahaye, B. Germain-Bonne, R{\'{e}}sultats n{\'{e}}gatifs en
ac\-c{\'{e}}\-l{\'{e}}\-ra\-tion de la convergence, Numer. Math., 35 (1980) 443-457.


\bibitem{rre}
R.P. Eddy, Extrapolation to the limit of a vector sequence, in {
Information Linkage between Applied Mathematics and Industry}, P.C.C. Wang
ed., Academic Press, New York, 1979, pp. 387-396.

\bibitem{eldendata}
L. Eld\'en, Numerical linear algebra in data mining, Acta Numer., 15 (2006) 327-84.

\bibitem{Evans}
C. Evans, S. Pollock, L. G. Rebholz, M. Xiao,  A Proof That Anderson Acceleration Improves the Convergence Rate in Linearly Converging
Fixed-Point Methods (But Not in Those Converging Quadratically), SIAM J. Numer. Anal., 58(1)  (2020) 788-810.

\bibitem{eyert96}
  V.~Eyert, A comparative study on methods for convergence
    acceleration of iterative vector sequences, J. Comput. Phys, 124 (1996) 271-285.

\bibitem{fang2009two}
H.R. Fang,  Y. Saad,  Two classes of multisecant methods for nonlinear acceleration, Numer. Linear Algebra Appl., 16(3) (2009) 197-221.

\bibitem{B}
A. Fu, J. Zhang, S. Boyd,  Anderson accelerated Douglas-Rachford splitting. SIAM J. Sci.
Comput., 42(6) (2020) A3560-A3583.

\bibitem{projproj}
A. Gal\'antai, {\it Projectors and Projection Methods}, Springer Science \& Business Media, 2003.

\bibitem{GSsnep} D.M. Gay, R.B. Schnabel, Solving systems of nonlinear equations by Broyden's method with projected updates,
in {\it  Nonlinear Programming}, Vol. 3,  O. Mangasarian, R. Meyer and S. Robinson, eds., Academic Press,  New York, 1978, pp. 245-281.


\bibitem{gek}
E. Gekeler, On the solution of systems of equations by the epsilon algorithm of Wynn, Math.
Comput., 26 (1972) 427-436.

\bibitem{GCVgolub}
G.H. Golub, M. Heath, G. Wahba, Generalized cross-validation as a method for choosing a good ridge parameter, Technometrics, 21(2) (1979)  215-23.

\bibitem{guSRRQR}
M. Gu, S.C. Eisenstat, Efficient algorithms for computing a strong rank-revealing QR factorization, SIAM J. Sci. Comput., 17(4) (1996) 848-869.


\bibitem{hjbbratu}
M. Hajipour, A. Jajarmi, D. Baleanu, On the accurate discretization of a highly nonlinear boundary value problem, Numer. Algorithms, 79(3) (2018)  679 - 95.

\bibitem{henrici}
P. Henrici, {\it Elements of Numerical Analysis}, Wiley, New York, 1964.


\bibitem{a4}
N.J. Higham, N. Strabi\'c, Anderson acceleration of the alternating projections
method for computing the nearest correlation matrix,
Numer. Algorithms, 72 (2016) 1021-1042.

\bibitem{rrrr}
K. Jbilou, H. Sadok, Some results about vector extrapolation methods and related fixed point
iteration, J. Comp. Appl. Math., 36 (1991) 385-398.

\bibitem{kelley2018numerical}
C.T. Kelley, Numerical methods for nonlinear equations, Acta Numerica, 27 (2018) 207-287.

\bibitem{hlf}
H. Le Ferrand, The quadratic convergence of the topological epsilon algorithm for systems of nonlinear equations,
Numer. Algorithms, 3 (1992) 273-284.

\bibitem{lupo2019convergence}
M. Lupo Pasini, Convergence analysis of Anderson-type acceleration of Richardson's iteration,  Numer. Linear Algebra Appl., 26(4) (2019) e2241.
\bibitem{rre1}
M. Me\`sina, Convergence acceleration for the iterative solution of $x=Ax+f$,
Comput. Methods Appl. Mech. Eng., 10 (1977) 165-173.

\bibitem{orte}
J.M. Ortega, W.C. Rheinboldt, {\it Iterative Solution of Nonlinear Equations in Several Variables}, Academic Press,
New York, 1970.

\bibitem{C}
W. Ouyang, J. Tao, A. Milzarek, B. Deng, Nonmonotone globalization for Anderson Acceleration using adaptive
regularization, arXiv:2006.02559 (2020).

\bibitem{scikit}
F. Pedregosa, G. Varoquaux, A. Gramfort, V. Michel, B. Thirion, O. Grisel, M. Blondel, P. Prettenhofer, R. Weiss,
V. Dubourg,  J. Vanderplas, Scikit-learn: machine learning in Python, J. Mach. Learn. Res., 12 (2011), 2825-30.


\bibitem{pollock1}
S. Pollock, L. G. Rebholz, M. Xiao, Anderson-accelerated convergence of Picard iterations for incompressible
Navier-Stokes equations, SIAM J. Numer. Anal., 57(2) (2019), 615-637.

\bibitem{puga}
B.P. Pugachev, Acceleration of convergence of iterative processes and a
method of solving systems of non-linear equations, USSR Comput. Maths.
Maths. Phys., 17 (5) (1978) 199-207.

\bibitem{pul1}
P. Pulay, Convergence acceleration in iterative sequences: the case of SCF iteration,
Chem. Phys. Lett., 73 (1980) 393-398.


\bibitem{RthesisQC} T. Rohwedder, {\it An analysis for some methods and algorithms of quantum chemistry}, PhD thesis, TU Berlin, 2010.

  \bibitem{RSDIIS} T. Rohwedder T, R. Schneider, An analysis for the DIIS acceleration method used in quantum chemistry calculations, J. Math. Chem., 49(9) (2011) 1889-1914.

\bibitem{aspre}
D. Scieur, A. D'Aspremont, F. Bach, Regularized nonlinear acceleration, Math. Program., 179 (2020) 47-83.


\bibitem{sha2}
D. Shanks, Non linear transformations of divergent and slowly convergent
sequences, J. Math. and Phys., 34 (1955) 1-42.

\bibitem{sidiprojection}
A. Sidi, Extrapolation vs. projection methods for linear systems of equations, J. Comput. Appl. Math., 22(1) (1988) 71-88.
\bibitem{sidisvd}
A. Sidi, SVD-MPE: An SVD-based vector extrapolation method of polynomial type, Applied Mathematics, 7 (2016) 1260-1278.
\bibitem{sidib}
A. Sidi, {\it Practical Extrapolation Methods. Theory and Applications},
Cambridge University Press, Cambridge, 2003.
\bibitem{sidibri}
A. Sidi, J. Bridger, Convergence and stability analyses for some vector extrapolation methods in the presence of
defective iteration matrices, J. Comp. Appl. Math., 22 (1988) 35–61.

\bibitem{sfs}
D.A. Smith, W.F. Ford, A. Sidi, Extrapolation methods for vector sequences, SIAM
Rev., 29 (1987) 199-233; Correction, SIAM Rev., 30 (1988) 623-624.
\bibitem{skel}
S. Skelboe, Computation of the periodic steady-state response to non linear networks by extrapolation
methods, IEEE Trans. Circuits Syst., 27 (1980) 161-175.

\bibitem{stef}
J.F. Steffensen, Remarks on iteration, Skand. Aktuarietidskr., 16 (1933) 64–-72.

\bibitem{toth2015convergence}
A. Toth,  C.T. Kelley, Convergence analysis for Anderson acceleration, SIAM J. Numer. Anal., 53(2) (2015) 805-819.
\bibitem{WalkerNi2011}
H.F. Walker, P. Ni, Anderson acceleration for fixed-point iterations,
SIAM J. Numer. Anal., 49, (2011) 1715-1735.
\bibitem{weni}
E.J. Weniger, Nonlinear sequence transformations for the acceleration of
convergence and the summation of divergent series, Comput. Phys. Rep., 10 (1989) 189-371.
\bibitem{wimp}
J. Wimp, {\it Sequence Transformations and their Applications},
Academic Press, New York, 1981.

\bibitem{zhangboyd}
J. Zhang,  B. O'Donoghue, S. Boyd, Globally convergent type-I Anderson acceleration for non-smooth fixed-point iterations,
arXiv:1808.03971 (2018).


\end{thebibliography}
\end{document}